\documentclass[11pt]{article}

\usepackage{amsfonts}
\usepackage{amsthm}
\usepackage{color}
\usepackage{graphicx}
 \usepackage{rotating}
\usepackage{url}

\newcommand\bA {\mathbf A}

\newcommand\bH {\mathbf H}

\newcommand\bO {\mathbf O}

\newcommand\bU {\mathbf U}
\newcommand\bV {\mathbf V}
\newcommand\bW {\mathbf W}

\newcommand\bY {\mathbf Y}
\newcommand\bZ {\mathbf Z}

\newcommand\wtbA {\widetilde{\bA}}

\newcommand\wtX {\widetilde{X}}
\newcommand\wtZ {\widetilde{Z}}

\newcommand\itF {{\mathcal{F}}}
\newcommand\itG {{\mathcal{G}}}
\newcommand\itH {{\mathcal{H}}}
\newcommand\itI {{\mathcal{I}}}

\newcommand\itU {{\mathcal{U}}}

\newcommand\bDelta {\mbox{\boldmath $\Delta$}}
\newcommand\bGa {\mbox{\boldmath $\Gamma$}}

\newcommand\bSi {\mbox{\boldmath $\Sigma$}}

\newcommand\bUpsi {\mbox{\boldmath $\Upsilon$}}

\newcommand\wgamma {\widehat{\gamma}}

\newcommand\wlam {\widehat{\lambda}}

\newcommand\wtheta {\widehat{\theta}}

\newcommand\wtau {\widehat{\tau}}

\newcommand\wbGa  {\widehat{\bGa}}

\newcommand\wbDelta {\widehat{\bDelta}}

\newcommand\wbSi {\widehat{\bSi}}

\newcommand\wbUps  {\widehat{\bUpsi}}

\newcommand\wtbGa {\widetilde{\bGa}}

\newcommand\wtbUps {\widetilde{\bUpsi}}

\def\real{\mathbb{R}}

\newcommand{\esp}{\mathbb{E}}
\newcommand{\prob}{\mathbb{P}}
\newcommand{\cov}{\mbox{\sc Cov}}
\newcommand{\var}{\mbox{\sc Var}}

\newcommand{\diag}{\mbox{\sc diag}}

\newcommand{\convprob  }{ \buildrel{p}\over\longrightarrow}

\newcommand{\convdist}{ \buildrel{D}\over\longrightarrow}

\newcommand{\trasp}{^{\mbox{\footnotesize \sc t}}}

\newcommand{\identidad}{\mbox{\bf I}}

\newcommand{\boot}{\mbox{\footnotesize\sc b}}
\newcommand{\weigh}{\mbox{\footnotesize\sc w}}
\newcommand{\perm}{\mbox{\footnotesize\sc p}}
\newcommand{\gaus}{\mbox{\footnotesize\sc g}}

\newcommand\noi{\noindent}

\def\square{\ifmmode\sqr\else{$\sqr$}\fi}
\def\sqr{\vcenter{
         \hrule height.1mm
         \hbox{\vrule width.1mm height2.2mm\kern2.18mm
\vrule width.1mm}
         \hrule height.1mm}}

\begin{document}


\title{Testing equality between several populations covariance operators}
\author{{ Graciela Boente},  {Daniela Rodriguez} {and} { Mariela Sued}\\
{\small \sl Facultad de Ciencias Exactas y Naturales, Universidad de Buenos Aires and
CONICET, Argentina} \\
{\small e--mail: gboente@dm.uba.ar\hskip1truecm drodrig@dm.uba.ar\hskip1truecm msued@dm.uba.ar}
}
\date{}
\maketitle

\begin{abstract}
In many situations, when dealing with several populations, equality of the covariance  operators is assumed. An important issue is to
study if this assumption holds before making other inferences. In this paper, we develop a test for comparing covariance operators of
several functional data samples. The proposed test is based on the Hilbert--Schmidt norm of  the difference between   estimated covariance operators. In particular, when dealing with two populations, the tests statistic is just the squared norm of the  difference between the two   covariance operators estimators. The    asymptotic behaviour of the test statistic under the null and under  local  alternatives  is obtained. Since the statistic null asymptotic distribution  does not allow to obtain easily its quantiles,   a bootstrap procedure to compute the critical values is considered. The performance of the test statistics for small sample sizes is illustrated through a Monte Carlo study.
 
%
\end{abstract}
  
 \section{Introduction}

 In many applications, we  study  phenomena that are continuous in time or space and can be considered as smooth curves or functions. Statistical procedures to deal with such functional data may be found, for instance, in Ramsay and Silverman (2005), Ferraty and  Vieu  (2006) and  Horv\'ath and Kokoszka (2012).
 On the other hand, when working with more than one population, as in the finite--dimensional case,  equality among the  covariance
 operators associated to each population is often assumed.  In the case of finite--dimensional data, tests for equality of covariance
 matrices have been extensively studied, see for example, Seber (1984) and Gupta and Xu (2006). This problem has been considered even for
 high dimensional data, i.e., when the sample size is smaller than the number of variables under study; we refer among others to Ledoit
 and Wolf  (2002) and Schott (2007).
 
 For functional data, most of the literature on hypothesis testing deals with tests on the mean function including the functional linear model, see, for instance, Fan and Lin (1998), Cardot \textsl{et al.} (2003),  Cuevas \textsl{et al.} (2004) and  Shen and Faraway (2004). However, as mentioned in Pigoli \textsl{et al.} (2014) the analysis of the covariance operator arises in many applied contexts. For instance,  Ferraty \textsl{et al.} (2007) considered tests for comparing groups of curves based on comparison of their covariances.  In this paper,  we rather focus on testing equality of the covariance operators of several functional samples. 
 By the Karhunen--Lo\'eve expansion, this is equivalent to testing if
 all the samples have the same set of functional principal components sharing also their size.  When considering only two populations, Benko \textsl{et al.} (2009), Panaretos \textsl{et al.} (2010) and Fremdt \textsl{et al.} (2013)  used this characterization to develop   test statistics.
 In particular, Benko \textsl{et al.} (2009) proposed two--sample bootstrap tests for specific aspects of the spectrum of functional data, such as the equality of a subset of eigenfunctions. On the other hand, Panaretos \textsl{et al.} (2010) and Fremdt \textsl{et al.} (2013) considered an approach based on the projection of the data over a suitable  chosen finite--dimensional space, such as that defined by the functional principal components.  The results in  Fremdt \textsl{et al.} (2013) generalized those provided in Panaretos \textsl{et al.}  (2010) which assume that  the processes  have a Gaussian distribution. More recently, Pigoli \textsl{et al.} (2014)  developed a two--sample test for comparing covariance operators using different   distances between covariance operators. Their procedure is based on a permutation test and assumes that the two samples have the same mean, otherwise, an approximate permutation test is considered  after the processes are centred using their sample means. A different approach was given by  Gaines \textsl{et al.} (2011), who defines an  univariate likelihood ratio test combined with Roy's union--intersection principle for testing the equality of two covariance operators and derives its asymptotic behaviour under the null and under a set of local alternatives converging to the null hypothesis at  a rate  $n^{1/2}$, where $n$ stands for the total sample size. 
 
 In this paper, we go one step further by studying the problem of more than two populations, that is,   we study if the covariance operators of $k$ independent samples are equal,  for $k\ge 2$.  Clearly, the permutation test  defined in Pigoli \textsl{et al.} (2014) can easily be adapted to the situation of more than two populations. However, up to our knowledge, the asymptotic behaviour of their test statistic under   local alternatives  has not been studied yet. One of the goals of this paper is not only to propose a test statistic to compare covariance operators of $k$ populations, but also to provide a theoretical framework which   clarifies the ability of the test statistic to detect root$-n$  local  alternatives and the rate of convergence of the detected alternatives.  Hence, our results extend the approaches based on distances between covariance operators estimators  given in the case of two independent samples to the several samples situation and provide a full asymptotic analysis not only under the null but also under local alternatives converging at a root$-n$ rate, which include, for instance, the functional common principal components model. 
 
 To fix ideas, we will begin by  describing the two sample situation and we will then generalize the test statistic to the situation in which $k>2$. Let
 us assume that we have two independent populations with covariance operators  $\bGa_1$ and  $\bGa_2$. Denote by $\wbGa_1$ and
 $\wbGa_2$ consistent estimators of $\bGa_1$ and  $\bGa_2$, respectively,  such as the sample covariance estimators studied in
 Dauxois \textsl{et al.} (1982). It is clear that under the standard null hypothesis $\bGa_1=\bGa_2$, the difference between the
 covariance  operator estimators should be small. For that reason, a test statistic based on the norm of
 $\wbGa_1-\wbGa_2$ may be helpful to study the hypothesis of equality. Different norms including the operator norm and the trace have been explored in  Pigoli \textsl{et al.} (2014) for the two sample case. Here, we will focus on the Hilbert--Schmidt norm which is a natural norm for covariance operators. In the general situation, one may consider the   norm of $\wbGa_1-\wbGa_j$, for $j=2,\dots, k$ to construct a test statistic.  It is worth noting that the asymptotic null distribution of the test statistic was   stated without proof in   Boente \textsl{et al.} (2011) corresponding to the peer--reviewed contribution to the International Workshop on Functional and Operatorial Statistics. This paper completes the results stated therein by generalizing the results to several populations and by providing the asymptotic distribution for root$-n$  local  alternatives as well as a numerical study conducted to illustrate the behaviour of the proposed test for finite samples.

 The paper is organized as follows. Section \ref{prel}  introduce the notation and review some basic concepts which are used in later
 sections. In Section \ref{stest}, we introduce  the  test statistics and derive its   asymptotic distribution under the null
 hypothesis.   An important issue is to describe a set of local alternatives that the proposed statistic is able to detect. For that purpose,   the asymptotic distribution under a set of local alternatives converging to the null hypothesis at rate $n^{1/2}$ which include  the functional common principal component model is studied in Section \ref{alt}.  A bootstrap calibration for the null distribution of the test statistic  is  described in Section \ref{bootstrap}. The results of a Monte Carlo study are summarized in Section \ref{simu}.
   Proofs are relegated to the Appendix.
 
 \section{Preliminaries and notation}\label{prel}
 From now on,    $\itH$  stands for a
 separable Hilbert space with inner
 product $\langle \cdot,\cdot\rangle$ and norm $\|u\|=\langle
 u,u\rangle^{1/2}$. Let $\bH:\itH\to \itH$ be a compact operator. 
 The operator $\bH$ is said to be a trace class operator if $\sum_{\ell=1}^\infty \langle \bH  u_\ell,   u_\ell\rangle<\infty$, for 
 $\{u_\ell :\ell\,\ge 1\}$ any orthonormal basis of $\itH$,
 while it is said to be Hilbert--Schmidt if  $\sum_{\ell=1}^\infty \|\bH u_\ell\|^2<\infty$.
 From now on, $\itF$ stands for the Hilbert space of  Hilbert--Schmidt operators over $\itH$ and
 $\bH^{*}$ will denote the adjoint of the operator $\bH$.  Given $\bH_1$,
 $\bH_2$ and $\bH$  Hilbert-Schmidt operators, the inner product in $\itF$ is defined as  $\langle \bH_1,
 \bH_2\rangle_{\itF} =\mbox{trace} (\bH_1^{*}\bH_2)
 =\sum_{\ell=1}^\infty \langle \bH_1 u_\ell, \bH_2 u_\ell\rangle$, while the
 norm equals $\|\bH\|_{\itF}=\langle \bH^{*} , \bH \rangle_{\itF}^{1/2}
 =\{\sum_{\ell=1}^\infty \|\bH u_\ell\|^2\}^{1/2}$, with $\{u_\ell :
 \ell\,\ge 1\}$ any orthonormal basis of $\itH$.
 These definitions are independent of the basis choice. Note that Hilbert--Schmidt operators
 have a countable number of
 eigenvalues, all of them being real when the operator is self--adjoint. Hence, choosing   as orthonormal basis
 the eigenfunctions  of $\bH$,  we get that for non--negative and self-adjoint operators $\| \bH
 \|^2_{\itF} =\sum_{\ell=1}^\infty \lambda^2_{\ell}$, where   $\{\lambda_{ \ell}: \ell\,\ge 1\}$ are
 the eigenvalues of $\bH$ ordered such
 that $\lambda_{ \ell}\ge \lambda_{ \ell+1}$.

 Let us consider independent  random elements $X_1,\dots,X_k$ in $\itH$ and assume that $\esp\|X_i\|^2<\infty$.  Denote by
 $\mu_i\in {\itH}$ the mean of $X_i$,  $\mu_i = \esp(X_i)$ and by
 $\bGa_i:{\itH}\to{\itH}$  the covariance operator of  $X_i$. Let
 $\otimes$ stand for the tensor product on $\itH$, e.g., for $u,
 v\in {\itH}$, the operator $u\otimes v:{\itH}\to \itH$ is
 defined as $(u\otimes v)w= \langle v,w\rangle u$. With this
 notation,  the covariance operator $\bGa_i$ can be written as
 $\bGa_i=\esp\{(X_i-\mu_i)\otimes(X_i-\mu_i)\}$.  The operator
 $\bGa_i$ is a linear, self--adjoint and compact operator with finite trace,
 so it is a Hilbert--Schmidt operator. From now on, we denote as $\{\phi_{i,\ell}: \ell\,\ge 1\}$
 the eigenfunctions of $\bGa_i$ related to the eigenvalues $\{\lambda_{i,\ell}: \ell\,\ge 1\}$,
 ordered as a non--increasing sequence, i.e., $\lambda_{i,\ell}\ge \lambda_{i,\ell+1}$.
 Note that the trace of $\bGa_i$ is given by $\sum_{\ell=1}^\infty \lambda_{i,\ell}$.
 
 When  ${\itH}=L^2({\itI})$  for some bounded interval $\itI$  and $\langle u,
 v\rangle=\int_{\itI} u(s)v(s) ds$, it is well known that the covariance operator is
 defined through the covariance function of $X_i$, $\gamma_i(s,t) =
 \cov(X_i (s), X_i(t))$, $s, t \in \itI$ as $(\bGa_i
 u)(t)=\int_{{\itI}} \gamma_i(s,t) u(s) ds$.  Besides, $\| \bGa_i
 \|^2_{\itF} =\sum_{\ell=1}^\infty \lambda^2_{i,\ell}=
  \|\gamma_i\|^2= \int_{\itI} \int_{\itI} \gamma_i^2(t,s) dt \, ds $.

 Our goal is to test whether the covariance operators $\bGa_i$  of
 several populations are equal or not. For that purpose, let us
 consider independent samples of each population, that is, let us
 assume that we have independent observations
 $X_{i,1},\cdots,X_{i,n_i}$, $1\le i \le k$ such that $X_{i,j}\sim X_i$, $1\le j\le n_i$.
 A natural way to estimate the covariance operators $\bGa_i$  is through their empirical versions. The sample
 covariance operator $\wbGa_i$ is defined as
 $$\wbGa_{i}=\frac{1}{n_i}\sum_{j=1}^{n_i}
 \left(X_{i,j}-\overline{X}_i\right)\otimes
 \left(X_{i,j}-\overline{X}_i\right)\;,$$
  where
 $\overline{X}_i=(1/n_{i})\sum_{j=1}^{n_i}X_{i,j}$. Dauxois \textsl{et
 al.} (1982) obtained the asymptotic behaviour of $\wbGa_i$. In
 particular, they have shown that, when $\esp(\|X_{i}\|^4)<\infty$,
  $\sqrt n_i\left({\wbGa_i}-\bGa_i\right)$ converges in distribution to
 a zero mean Gaussian random element of $\cal F$, $\bU_i$, with
 covariance operator $\bUpsi_i$  given by
 \begin{eqnarray}
 \bUpsi_i&=&\sum_{m,r,o,p} \lambda_{i m}^{1/2}\lambda_{i r}^{1/2} \lambda_{i o}^{1/2} \lambda_{i p}^{1/2}
 \esp\left(f_{im}f_{ir}f_{io}f_{ip}\right)\,\phi_{i,m}\otimes\phi_{i,r}
 \tilde{\otimes} \phi_{i,o}\otimes\phi_{i,p}\nonumber\\
 & -&\sum_{m,r}
 \lambda_{im}\lambda_{ir}\,\phi_{i,m}\otimes\phi_{i,m}
 \tilde{\otimes} \phi_{i,r}\otimes\phi_{i,r} \label{varasintgamai}
 \end{eqnarray}
 where $\tilde{\otimes}$ stands for the tensor product in $\itF$
 and, as mentioned above, $\{\phi_{i,\ell}: \ell\,\ge 1\}$ is an
 orthonormal basis of eigenfunctions of $\bGa_i$ with associated
 eigenvalues  $\{\lambda_{i,\ell}: \ell\,\ge 1\}$  such that
 $\lambda_{i,\ell}\ge\lambda_{i,\ell+1}$. The random variables $f_{im}$ are the
 standardized coordinates of $X_i-\mu_i$ on the basis
 $\{\phi_{i,\ell}:\ell\,\ge 1\}$, that is, $f_{im}=\lambda_{i,m}^{-1/2}\,\langle
 X_i-\mu_i,\phi_{i,m} \rangle $. Note that
 $\esp(f_{im})=0$. Using that $\cov\left(\langle u, X_i-\mu_i\rangle,
 \langle v, X_i-\mu_i\rangle\right)=\langle u, \bGa_i v\rangle$, we
 get that  $\esp(f_{im}^2)=1$, $\esp(f_{im}\;f_{is})=0$ for $m\ne s$.
 In particular, the Karhunen-Lo\'eve expansion leads to
 \begin{equation}
 X_i\;=\mu_i+\;\sum_{\ell=1}^\infty \lambda_{i,\ell}^{\frac
 12}\,f_{i\ell}\,\phi_{i,\ell}\;.
 \label{kl}
 \end{equation}
 It is worth noticing that, since $\esp \|\bU_i\|_\itF^2<\infty$,
 $\bUpsi_i$ is a linear  operator over $\cal F$ with finite trace, so it is also a Hilbert--Schmidt operator. Thus, any linear combination  of the operators $\bUpsi_i$,
 $\bUpsi=\sum_{i=1}^k a_i \bUpsi_i$, with $a_i\ge 0$, will be trace class operator, that is, if $\{\theta_\ell \}_{\ell\,\ge
 1}$ stand for the eigenvalues of $\bUpsi$ ordered in decreasing
 order, we have that $\theta_\ell \ge 0$ and $\sum_{\ell\,\ge 1} \theta_\ell
 <\infty$. This property will be used later in Theorem \ref{stest}.1.
 
 When ${\itH}=L^2({\itI})$, smooth estimators, ${\wbGa_{i,h}}$, of the covariance operators were
 studied in Boente and Fraiman (2000).  The smoothed operator is the operator induced by the smooth covariance function
 $$\wgamma_{i,h}(t,s)=\frac{1}{n_1}\sum_{j=1}^{n_i} \left(X_{i,j,h}(t)-\overline{X}_{i,h}(t)\right)
 \left(X_{i,j,h}(s)-\overline{X}_{i,h} (s)\right)\;,$$ where
 $X_{i,j,h}(t)=\int_{\itI} K_h(t-x)X_{i,j}(t) dt $ are the smoothed
 trajectories, $K_h(\cdot)=h^{-1}K(\cdot/h)$ is a nonnegative kernel
 function, and $h$ a smoothing parameter.  Boente and Fraiman (2000)
 have shown that the smooth estimators have
 the same asymptotic distribution as the empirical version, under mild conditions.

 \section{The test statistic}\label{stest}
 
 To motivate our test statistic, we first consider the two sample setting, that is, the problem of testing the hypothesis
 \begin{equation}
 H_0: \bGa_1=\bGa_2\quad\mbox{against}\quad H_1: \bGa_1\not=\bGa_2\;\;,
 \label{test}
 \end{equation}
 from two independent samples $X_{1,1},\cdots,X_{1,n_1}$ and $X_{2,1},\cdots,X_{2,n_2}$.
 A natural approach is to consider $\wbGa_i$ as the empirical
 covariance operators of each population  and construct  a statistic $T_n$ based on the difference between  the covariance operators
 estimators, i.e.,  to define
 \begin{equation}
 T_n=n \|\wbGa_1 -\wbGa_2 \|_\itF^2\, ,
 \label{Tnk=2}
 \end{equation}
  where $n=n_1+n_2$, ${n_i}/n\to \tau_i$ with $\tau_i\in (0,1)$. As mentioned in  Pigoli \textsl{et al.} (2014), the null hypothesis can be written as $d(\bGa_1, \bGa_2)=\|\bGa_1 -\bGa_2 \|_\itF=0$ while the alternative corresponds to $\|\bGa_1 -\bGa_2 \|_\itF>0$. Thus, if  $\wbGa_j$ are consistent estimators of $\bGa_j$ for $j=1,2$, any test  based on the distance $d(\wbGa_1, \wbGa_2)$ between will be consistent.

 To generalize the procedure to several populations, let $\bGa_i$ stand for  the covariance operator  of the $i-$th population. We wish to test the null hypothesis
 \begin{equation}
 H_0: \bGa_1=\dots=\bGa_k\quad\mbox{against}\quad H_1: \exists \; i\neq j \mbox{ such that }\bGa_i\not=\bGa_j\;.
 \label{testk}
 \end{equation}
 Note that the null hypothesis is equivalent to $\sum_{j=2}^k\|{\bGa}_j-{\bGa}_1\|_\itF^2=0$, which allows to 
 construct a consistent test using consistent covariance operator estimators. To be more precise, 
 let $X_{i,1},\cdots,X_{i,n_i}$, $1\le i \le k$, be independent samples, $n=n_1+\dots+n_k$ and assume 
 that $n_i/n\to \tau_i$, $0<\tau_i<1$, $\sum_{i=1}^k \tau_i=1$. Denote with   $\wbGa_i$   the sample covariance operator  of $i-$th population.   A natural generalization of the statistic defined in (\ref{Tnk=2}) is to consider the  test statistic
 \begin{equation}
 T_{k,n}=n \sum_{j=2}^k\|{\wbGa}_j-{\wbGa}_1\|_\itF^2\,.
 \label{tkn}
 \end{equation}
 To define the test we need the asymptotic distribution of $T_{k,n}$ under the null hypothesis, which is derived in Corollary \ref{stest}.1.

 The following result allows  to study, under the null hypothesis, the asymptotic behaviour of $n \sum_{j=2}^k\|{\wtbGa}_j-{\wtbGa}_1\|_\itF^2$ when considering a general class of covariance estimators $\wtbGa_i$ rather than the sample covariance operators.

 \vskip0.1in
 \noi \textbf{Theorem \ref{stest}.1.} \textsl{Let $X_{i,1} ,\cdots,X_{i,n_i} $, for $1\le i\le k$, be
 independent observations from $k$ independent
 distributions in $\itH$, with mean $\mu_i$ and
 covariance operator $\bGa_i$.  Assume that ${n_i}/n\to \tau_i$ with $\tau_i\in (0,1)$
 where $n=\sum_{i=1}^k n_i$. Let $\wtbGa_i$ be the independent   estimators of  the $i-$th population covariance
 operator such that $\sqrt n_i\left({\wtbGa_i}-\bGa_i\right)\convdist
 \bU_i$,  with $\bU_i$ a zero mean Gaussian random element with
 covariance operator $\bUpsi_i$.
 Denote $ \bUpsi_{\weigh}=\; \left(\bUpsi_{\weigh, 1},\dots,\bUpsi_{\weigh, k-1}\right)$ the trace
 operator $ \bUpsi_{\weigh}:{\itF}^{k-1}\to {\itF}^{k-1}$ with $i-$th component defined
 as
 \begin{equation}
 \bUpsi_{\weigh, i}(u_1,\dots,u_{k-1})= \frac{1}{\tau_{i+1}}\bUpsi_{i+1}({u_i})+\frac{1}{\tau_1}\bUpsi_1 \left(\sum_{\ell=1}^{k-1}u_\ell\right)\qquad \qquad \mbox{for } 1\le i\le k-1\,.
 \label{upsiweight}
 \end{equation}
 Let
 $\{\theta_\ell\}_{\ell\,\ge 1} $ stand for the sequence of
 eigenvalues of $  \bUpsi_{\weigh}$ ordered in decreasing order. Then, we have that
 $$
 n \sum_{j=2}^k\|(\wtbGa_j-\bGa_j)-(\wtbGa_1-\bGa_1)\|_\itF^2\buildrel{{\cal D}}\over\longrightarrow  \sum_{\ell\,\geq 1}\theta_\ell  Z_\ell^2\,,
   $$
  with $Z_\ell\sim N(0,1)$  independent. In particular, if $H_0: \bGa_1=\dots=\bGa_k$ holds, we have that
 $n \sum_{j=2}^k\| \wtbGa_j -\wtbGa_1 \|_\itF^2\buildrel{{\cal D}}\over\longrightarrow  \sum_{\ell\,\geq 1}\theta_\ell  Z_\ell^2$.
 }
 
 \vskip0.1in
 
 When  $\esp(\|X_{i}\|^4)<\infty$, the results in Theorem  \ref{stest}.1 apply in particular to the sample covariance operator, i.e., when $\wtbGa_i=\wbGa_i$, leading to the asymptotic  distribution of $T_{k,n}$ under the null hypothesis stated in Corollary \ref{stest}.1. However, it also allows to use other  covariance estimators to define the test statistic, such as the smooth ones  ${\wbGa_{i,h}}$ defined in Boente and Fraiman (2000).
 
 \vskip0.1in
 \noi \textbf{Corollary \ref{stest}.1.} \textsl{Let $X_{i,1} ,\cdots,X_{i,n_i} $, for $1\le i\le k$, be
 independent observations from $k$ independent
 distributions in $\itH$, with mean $\mu_i$ and
 covariance operator $\bGa_i$ such that $\esp(\|X_{i}\|^4)<\infty$.
 Let $\wbGa_i$ be the sample covariance operator of the $i-$th
 population. Assume that ${n_i}/n\to \tau_i$ with $\tau_i\in (0,1)$
 where $n=\sum_{i=1}^k n_i$. Denote $ \bUpsi_{\weigh}=\; \left(\bUpsi_{\weigh, 1},\dots,\bUpsi_{\weigh, k-1}\right)$ the trace
 operator $ \bUpsi_{\weigh}:{\itF}^{k-1}\to {\itF}^{k-1}$ where  $\bUpsi_{\weigh, i}$ is defined in (\ref{upsiweight}) with 
 $\bUpsi_{i}$  given in (\ref{varasintgamai}). Let
 $\{\theta_\ell\}_{\ell\,\ge 1} $ stand for the sequence of
 eigenvalues of $  \bUpsi_{\weigh}$ ordered in decreasing order. Under $H_0: \bGa_1=\dots=\bGa_k$, we have
 \begin{equation}
 n \sum_{j=2}^k\|\wbGa_j-\wbGa_1\|_\itF^2\buildrel{{\cal D}}\over\longrightarrow  \sum_{\ell\,\geq 1}\theta_\ell  Z_\ell^2\,,
 \label{distkpob}
 \end{equation}
  with $Z_\ell\sim N(0,1)$  independent.
 }

 \vskip0.1in
 \noi \textbf{Remark \ref{stest}.1.}
 \begin{itemize}
 \item[a)] Note that the fact that $\esp(\|X_{i}\|^4)<\infty$ entails that $\esp(\|(X_{i}-\mu_i)\otimes (X_{i}-\mu_i) \|^2)<\infty$, 
  so  $\bUpsi_{i}$, the covariance operator of $(X_{i}-\mu_i)\otimes (X_{i}-\mu_i)$, is well defined and
  $\sum_{\ell\ge 1}\theta_{\ell}<\infty$. Hence, for any $q_n$   a sequence of integers such that $q_n\to \infty$, the
 sequence ${\itU}_n=\sum_{\ell= 1}^{q_n} \theta_\ell  Z^2_i$ is Cauchy
 in $L^2(\prob)$, so the limit  ${\itU}=\sum_{\ell\,\geq 1} \theta_\ell  Z^2_\ell$
 is well defined. In fact, analogous arguments to those considered
 in Neuhaus (1980) allow to show that the series converges almost
 surely. Moreover, since $Z^2_1\sim \chi_1^2$, ${\itU}$ has a continuous
 distribution function $F_{\itU}$ which entails that  the distribution
 function of ${\itU}_n$,  $F_{{\itU}_n}$, converges to  $F_{\itU}$ uniformly (see, for instance, shown in Lemma 2.11 in  Van der Vaart, 2000).
 
 \item[b)] It is worth noticing  that Corollary \ref{stest}.1 is a natural extension of its analogous in the finite--dimensional case. To be more precisely,  let $\bZ_{ij}\in \real^p$ with $1\leq i\leq k$ and  $1\leq j\leq n_i$ be independent random  vectors and let  $\wbSi_i$ be their  sample  covariance matrix. Then,  $\sqrt{n_i} \bV_i=\sqrt{n_i} (\wbSi_i-\bSi_i)$ converges to a multivariate normal distribution with mean zero and covariance matrix $\Upsilon_i$.
 Let
 $$\bA=\left(
   \begin{array}{ccccc}
     -\identidad_{p } & \identidad_{p} & 0 & \dots & 0 \\
     -\identidad_{p }& 0 &\identidad_{p} & \dots & 0 \\
     \vdots & \vdots & \vdots & \vdots & \vdots \\
     -\identidad_{p }&  0 & \dots & 0& \identidad_{p} \\
   \end{array}
 \right)$$
 where $\identidad_{p }$ stands for the   identity matrix of order $p$. Then, straightforward calculations allow to show that  $\sqrt{n} \bA(\bV_1,\dots,\bV_k)\trasp\convdist N(0,\Upsilon)$ where
 $$\Upsilon=\left(
   \begin{array}{cccc}
     {\tau_1}^{-1}\Upsilon_1+ {\tau_2}^{-1}\Upsilon_2& {\tau_1}^{-1}\Upsilon_1& \dots & {\tau_1}^{-1}\Upsilon_1\\
     {\tau_1}^{-1}\Upsilon_1& {\tau_1}^{-1}\Upsilon_1+ {\tau_3}^{-1}\Upsilon_3 & \dots & {\tau_1}^{-1}\Upsilon_1\\
    \vdots & \vdots & \vdots & \vdots \\
    {\tau_1}^{-1}\Upsilon_1& {\tau_1}^{-1}\Upsilon_1 & \dots & {\tau_1}^{-1}\Upsilon_1+ {\tau_k}^{-1}\Upsilon_k\\
   \end{array}
 \right)$$ Therefore, under the null hypothesis of equality of the
 covariance matrices $\bSi_i$, we have that $
 n\sum_{i=2}^k\|\wbSi_i-\wbSi_1\|^2=\|\sqrt n \bA\bV\|^2\convdist
 \sum_{\ell=1}^{k p^4} \theta_\ell  Z_\ell^2$ where $\bV=(\bV_1,\dots,\bV_k)$ and
 $\theta_1,\theta_2,\dots, \theta_{k p^4}$ are the eigenvalues of
 $\Upsilon$. Note that the matrix $\Upsilon$ is the finite
 dimensional version of the covariance operator $ \bUpsi_{\weigh}$.
 \item[c)] From
 Corollary \ref{stest}.1 we have that, under the null hypothesis $H_0: \bGa_1=\dots=\bGa_k$, the test statistic $T_{k,n}=n \sum_{j=2}^k\|{\wbGa}_j-{\wbGa}_1\|_\itF^2,\convdist  {\itU}=\sum_{\ell\,\geq 1} \theta_\ell  Z^2_\ell$. Hence, an asymptotic test  may be based on $T_{k, n}$ rejecting for large values of
 $T_{k, n}$.  To obtain the critical values, the
 distribution of ${\itU}$ and thus, the eigenvalues of $ \bUpsi_{\weigh}$
  need to be estimated. In particular, when $k=2$, the test statistic $T_{k,n}$ equals $T_n=n \|\wbGa_1 -\wbGa_2 \|_\itF^2$ and $ \bUpsi_{\weigh}=\tau_1^{-1}\bUpsi_1+\tau_2^{-1}\bUpsi_2$. As
 mentioned above, the distribution function of ${\itU}$
 can be uniformly approximated by that of ${\itU}_n$ and so, the
 critical values can be approximated by the $(1-\alpha)-$quantile
 of  ${\itU}_n$. Gupta and Xu (2006) provide an approximation for
 the distribution function of any finite mixture of $\chi_1^2$
 independent random variables that can be used in the computation of
 the $(1-\alpha)-$quantile of  $\sum_{\ell= 1} ^{q_n}
 \wtheta_\ell Z^2_\ell$, where $\wtheta_\ell$ are
 estimators of $\theta_\ell $. It is also worth noticing that,    under
 $H_0: \bGa_1=\dots=\bGa_k$, the operator $\bUpsi_i$ given in
 (\ref{varasintgamai}) reduces to
 $$\bUpsi_i\!=\!\!\sum_{m,r,o,p} \lambda_m^{1/2} \lambda_r^{1/2}\lambda_o^{1/2}\lambda_p^{1/2} \esp[f_{im}f_{ir}f_{io}f_{ip}]\,\phi_{m}\otimes\phi_{r} \tilde{\otimes}
 \phi_{ o}\otimes\phi_{p} \!-\!\sum_{m,r} \lambda_{ m}\lambda_{r}\,\phi_{ m}\otimes\phi_{ m} \tilde{\otimes} \phi_{r}\otimes\phi_{r}\, $$
  for $i=1,\dots, k$, where, for the sake of simplicity, we denote as $\lambda_m$   the $m-$th largest eigenvalue of $\bGa_1$ and $\phi_m$ its corresponding eigenfunction.
  \newline
 In particular, if all the populations have the same underlying distribution except for the mean and covariance operator, as it happens when comparing the covariance operators of Gaussian processes, the random functions $f_{im}$, $i=2,\dots, k$, have the same distribution as $f_{1m}$,  so, in this case, $\bUpsi_1= \bUpsi_i$, for $i=2,\dots, k$, under $H_0$.
 
 \item[d)] Assume that the processes are  Gaussian, then using that $\esp (f_{im}f_{ir}f_{io}f_{ip} )$ equals 1 when pairs of indices are equal, 3 when $m=r=o=p$ and 0 otherwise, we have that, under the null hypothesis
 \begin{eqnarray*}
 \bUpsi_i=\bUpsi_1&=&\sum_{i\not=j} \lambda_i \lambda_j\,\phi_{i}\otimes\phi_{j} \tilde{\otimes}\phi_{j}\otimes\phi_{i}
 +\sum_{i\not=j} \lambda_i \lambda_j\,\phi_{i}\otimes\phi_{j} \tilde{\otimes}\phi_{i}\otimes\phi_{j}
  + 2\sum_{i} \lambda_i^2 \,\phi_{i}\otimes\phi_{i} \tilde{\otimes}\phi_{i}\otimes\phi_{i}\\
  &=& 2\sum_{i} \lambda_i^2 \,\phi_{i}\otimes\phi_{i} \tilde{\otimes}\phi_{i}\otimes\phi_{i} +\sum_{i<j} \lambda_i \lambda_j (\phi_i\otimes\phi_j+\phi_j\otimes\phi_i)\tilde{\otimes}(\phi_i\otimes\phi_j+\phi_j\otimes\phi_i)\,.
 \end{eqnarray*}
 Using that  $\phi_i\otimes\phi_i$ and
  $(\phi_i\otimes\phi_j+\phi_j\otimes\phi_i)/\sqrt{2}$, for  $i<j$,
  constitutes a complete orthonormal basis of the space of self--adjoint Hilbert--Schmidt operators, we conclude that they are the eigenfunctions of
  $\bUpsi_1$ associated to the eigenvalues  $2\lambda_i^2$ and $2\lambda_i\lambda_j$,
  respectively.
 Furthermore, if   $\tau_i=1/k$ for $i=1\dots, k$, we get that   for $1\le i\le k-1$
   \begin{equation}
 \bUpsi_{\weigh, i}(u_1,\dots,u_{k-1})= k\left[\bUpsi_1(u_i)+\bUpsi_1 \left(\sum_{\ell=1}^{k-1}u_\ell\right)\right]= k\left[\bUpsi_1(u_i)+\sum_{\ell=1}^{k-1}\bUpsi_1 \left(u_\ell\right)\right]\,,
 \label{upsiweight2}
 \end{equation}
 which entails that  $\theta_{i,i}=2k^2\lambda_i^2$ and $\theta_{i,j}=2k^2\lambda_i\lambda_j$, for $i<j$,
 are eigenvalues of $\bUpsi_{\weigh}=(\bUpsi_{\weigh, 1},\dots,\bUpsi_{\weigh, k-1})$, related to the eigenfunctions
 $v_{i,i}=(\phi_i\otimes\phi_i, \dots, \phi_i\otimes\phi_i)$ and
 $v_{i,j}=((\phi_i\otimes\phi_j+\phi_j\otimes\phi_i)/\sqrt{2},\dots,(\phi_i\otimes\phi_j+\phi_j\otimes\phi_i)/\sqrt{2})$,
 respectively. On the other hand, if $\alpha$ is an eigenvalue of  $\bUpsi_{\weigh}$,   $\alpha/k^2$
 is an eigenvalue of $\bUpsi_1$, meaning that we have obtained all the eigenvalues of $\bUpsi_{\weigh}$.
 \end{itemize}

 \section{Behaviour under local alternatives}{\label{alt}}
 
 In this section,  we study the behaviour of the test statistic $T_{k,n}$
 under a set of local alternatives. It is clear that, as in the multivariate situation, there are many ways in which
 the covariance operators may differ, one of them being the functional common principal model in which discrepancies from the
 null hypothesis arise only in the eigenvalues and not in the eigenfunctions of the covariance operators. Our results include that setting but also a situation in which the processes can be written as   sums of two independent processes, one of them having the same covariance operator along populations.
  
 We decided to keep fixed  the distribution of the first population, while that of the remaining ones will depend on the sample size,
 in such a way that for each fixed $n$ the alternative assumption holds but, as is usual for local alternatives, when the sample sizes increase, the alternatives considered converge to the null hypothesis at a given rate. 
 To avoid burden notation, when it is clear, in this section  we will use   $X_{i,j}$  to denote the observations under the local alternatives   $X_{i,j}^{(n)}$,  for  $1\leq j\leq n_i$, $2\leq i\leq k$. Similarly, we denote as $ X_i$, instead of $X_i^{(n)}$, the random element with common distribution, that is,   $X_{i,j}\sim X_i$.

  As in Section \ref{stest}, the following result present a general framework which allows to study the distribution of the test  statistic under root$-n$ local alternatives.   Theorem \ref{alt}.1 together with Propositions \ref{alt}.1a) and \ref{alt}.2a) allows  to derive the behaviour of the test statistic $T_{k,n}$ under the local alternatives described above. However, Theorem \ref{alt}.1 may also be applied when considering covariance estimators other than the sample covariance estimators. 
 
 \vskip0.1in
 \noi \textbf{Theorem \ref{alt}.1.} \textsl{Let $X_{i,1} ,\cdots,X_{i,n_i} $ for $i=1,\dots, k$ be independent observations
 from $k$ independent distributions  in $\itH$,  with  covariance
 operator $\bGa_i$ such that, for $i\ge 2$, $\bGa_i=\bGa_{i,n}=\bGa_1+n^{-1/2}\bDelta_i$. Assume that $\bDelta_i$ is a  self--adjoint trace operator such that $\bGa_{i,n}$ is non--negative. Denote as $\bDelta^{(k-1)}=(\bDelta_2,\dots, \bDelta_k)\trasp\in \itF^{k-1}$, $n=\sum_{i=1}^k n_i$ and
 assume    that ${n_i}/n\to \tau_i\in (0,1)$. Let
 $\wtbGa_i$   be the independent   estimators of  the $i-$th population covariance
 operator such that, for $1\le i\le k$, $\sqrt n_i\left({\wtbGa_i}-\bGa_1\right)\convdist \bU_i+ \tau_i^{1/2}\,\bDelta_i$ where $\bU_i$ is a zero mean Gaussian random element with covariance operator $\bUpsi_i$ and $\bDelta_1=\bO$ stands for the null operator.
   Define $\bUpsi_{\weigh}=(\bUpsi_{\weigh, 1}, \dots, \bUpsi_{\weigh, k-1})$ where  $\bUpsi_{\weigh, i}$ is given in (\ref{upsiweight}) and let  $\{\upsilon_\ell\}_{\ell\ge 1}$ be an orthonormal basis of eigenfunctions of $\bUpsi_{\weigh}$ related to the eigenvalues $\{\theta_\ell\}_{\ell \ge 1}$ ordered in decreasing order.   Then,
 $$n \sum_{i=2}^k\|\wtbGa_i -\wtbGa_1 \|_\itF^2\convdist  \sum_{\ell\,\geq 1} \theta_\ell  \left(Z_\ell+\frac{\eta_\ell}{\sqrt{\theta_\ell} }\right)^2\;,$$
  where $Z_\ell$ are independent and $Z_\ell \sim N(0,1)$ and
 $\eta_\ell=\langle \bDelta^{(k-1)},v_\ell\rangle_{\itF^{k-1}}$, i.e., $\bDelta^{(k-1)}=\sum_{\ell\,\ge 1}\eta_\ell \upsilon_\ell$.
 }
 \vskip0.1in
 
 Note that the requirement that $\bDelta_i$ is a  self--adjoint trace operator is needed to guarantee that $\bGa_i=\bGa_{i,n}$ is a valid covariance operator. Besides, since $\bDelta_i$ has finite trace, we have that $ \bDelta^{(k-1)}\in \itF^{k-1}$, so $\sum_{\ell\,\ge 1} \eta_\ell^2<\infty$.

   \vskip0.1in
 
  As mentioned at the beginning of this section, we will consider  two scenarios  where the assumptions of Theorem \ref{alt}.1 are satisfied. 
 The first one is a generalization of Example 2.2 in Gaines \textsl{et al.} (2011)  and assumes that, for $i=2,\dots,k$, the observations from the $i$-th population 
 can be written as the sum of two independent processes, the first one  having the same covariance operator as $X_1$.  Namely, we assume that 
 \begin{equation}
 \label{contaminado}
  X_{i,j}\;=X_{i,j}^{(n)}=W_{i,j}+n^{-1/4} R_{i,j}\;,\quad \mbox{ for } 2\leq i\leq k, 
 \end{equation}
 where $W_{i,j}$, $R_{i,j}$   are independent and such that $W_{i,j}\sim W_i$, $R_{i,j}\sim R_i$  and $W_{i}$ has the same covariance operator as $X_1$, for $1\le i\le k$.
  Notice that  the distribution of the  term $R_{i}$   is free to   vary across populations, for $2\leq i\leq k$, as well as the distribution 
 of $W_i$ as far as  $W_i$ and $X_1$   share the same covariance operator.  
 
  From now on, let  $\{\phi_\ell\}_{\ell\ge 1}$ be the eigenfunctions of $\bGa_1$,  the covariance operator of $X_1$, and   denote $\lambda_\ell$  the  eigenvalues of $\bGa_1$ related to  $\phi_\ell$, that is, we   omit  the subscript $1$ in $\lambda_{1,\ell}$  and $\phi_{1,\ell}$. 
 
 \vskip0.1in
 \noi  \textbf{Proposition \ref{alt}.1.} \textsl{Let $X_{i,1} ,\cdots,X_{i,n_i} $, $i=1,\dots, k$ be independent observations 
 from $k$ independent distributions  in $\itH$ such that (\ref{contaminado}) holds. Assume  that  ${n_i}/n\to \tau_i\in (0,1)$ with $n=\sum_{i=1}^k n_i$,
  $\esp(\|X_{1}\|^4)<\infty$  and that,  for $2\leq i\leq k$,   $\esp(\|W_{i}\|^4)<\infty$  and $\esp(\|R_{i}\|^4)<\infty$. 
 Let $\bDelta_i$ be the covariance operator   of 
  $R_{i}$,  for $i=2,\dots, k$ and assume that $\bGa_1=\esp\{(X_1-\mu_1)\otimes(X_1-\mu_1)\}$ is also the  covariance operator of $W_{i}$, for $i=2,\dots,k$. Denote as $\wbGa_i$ the sample covariance operator of  the $i-$th population.  
 Then, we have that $\sqrt n_i\left({\wbGa_i}-\bGa_1\right)\convdist \bU_i+ \tau_i^{1/2}\,\bDelta_i$ with $\bU_i$ a
  zero mean Gaussian random element with
 covariance operator $\bUpsi_i$ given in (\ref{varasintgamai}), that is,
 \begin{equation}
 \bUpsi_i=\!\!\sum_{m,r,o,p} \lambda_{m}^{1/2}\lambda_{r}^{1/2}\lambda_{o}^{1/2}\lambda_{p}^{1/2}  
 \esp[f_{im}f_{ir}f_{io}f_{ip}]\,\phi_{m}\otimes\phi_{r} \tilde{\otimes}
 \phi_{o}\otimes\phi_{p}-\sum_{m,r} \lambda_{m}\lambda_{r}\,\phi_{ m}\otimes\phi_{ m} \tilde{\otimes}
 \phi_{r}\otimes\phi_{r}\,, 
 \label{UPSIi1}
 \end{equation}
  where $f_{im}$ are the standardized coordinates of $W_i-\esp(W_i)$ on the basis $\{\phi_{\ell}:\ell\,\ge 1\}$, i.e., $\lambda_{\ell}^{\frac 12}\, f_{i\ell}=\langle W_i-\esp(W_i), \phi_\ell\rangle$.
   }

 \vskip0.1in
  Note that  if $W_{i}$ has the same distribution as $X_1$, we have that $\bUpsi_i=\bUpsi_1$.

 \vskip0.2in
  The second model  for local alternatives to be considered in this section is the functional common principal model.
 These alternatives include, as a particular case, 
 alternatives following the proportional model $\bGa_{i,n}=(1+\rho_i/\sqrt{n}) \bGa_1$.
 For details on the functional principal
 component model we refer to Benko \textsl{et al.} (2009) and
 Boente \textsl{et al.} (2010), for instance.
 By assuming   
 local alternatives satisfying  a functional  common principal model, we get that  the processes  $X_{i}$, $1\le i\le k$, can be
 written as
 \begin{equation}
 X_{1}\;=\mu_1+\;\sum_{\ell=1}^\infty \lambda_{\ell}^{\frac 12}\,f_{1\ell}\,\phi_\ell \quad \mbox{and} \quad X_{i}\;=X_i^{(n)}=\mu_i+\;\sum_{\ell=1}^\infty {\lambda_{i,\ell}^{(n)}}^{\frac 12}\,f_{i\ell}\,\phi_\ell, \mbox{ for } i\ge 2
 \label{fcpc}
 \end{equation}
 with $\lambda_{1}\ge \lambda_{2}\ge \dots\ge 0$,
 ${\lambda_{i,\ell}^{(n)}}\to \lambda_{\ell}$ at a given rate, while
 $f_{i\ell}$ are random variables such that $\esp(f_{i\ell})=0$,
 $\esp(f_{i\ell}^2)=1$, $\esp(f_{i\ell}\;f_{is})=0$ for $\ell\ne s$. 
 
  Proposition \ref{alt}.2 gives the asymptotic behaviour of the sample covariance operators when choosing  $\lambda_{i,\ell}^{(n)}=\lambda_{\ell} (1+n^{-1/2}\Delta_{i,\ell})$   in  (\ref{fcpc}). It is worth noting that,  if  $  (1+n^{-1/2}\Delta_{i,\ell})\ge 0$ and  some additional  conditions on $\Delta_{i,\ell}$ to be stated below are fulfilled, then $\bGa_i=\bGa_{i,n}=\bGa_1+n^{-1/2}\bDelta_i$, for $i\ge 2$,
 where $\bDelta_i=\sum_{\ell\,\ge 1}\Delta_{i,\ell} \lambda_{\ell} \phi_\ell\otimes \phi_\ell$. Hence, Proposition \ref{alt}.2 together with Theorem \ref{alt}.1 lead to the asymptotic behaviour of the test statistic $T_{k,n}$ under local alternatives following a functional common principal model, as stated in Corollary \ref{alt}.1.

 \vskip0.2in
 \noi \textbf{Proposition \ref{alt}.2.} \textsl{Let $X_{i,1} ,\cdots,X_{i,n_i} $, $i=1,\dots, k$, be independent observations 
 from $k$ independent distributions  in $\itH$, such that $X_{i,j}\sim X_i$.  Assume that $X_i$ satisfy (\ref{fcpc}) with
 $\lambda_{i,\ell}^{(n)}=\lambda_{\ell} (1+n^{-1/2}\Delta_{i,\ell})$ and that 
  ${n_i}/n\to \tau_i\in (0,1)$ where   $n=\sum_{i=1}^k n_i$. 
 Let $\wbGa_i$
 be the sample covariance operator of  the $i-$th population.
 Furthermore, assume that $\esp(\|X_{1}\|^4)<\infty$, $\sigma^2_{4,i,\ell}=\esp(f^4_{i\ell})<\infty$,
 $\sum_{\ell=1}^\infty \lambda_\ell |\Delta_{i,\ell}| <\infty$, $\sum_{\ell=1}^\infty \lambda_\ell \Delta_{i,\ell}^2
 \sigma_{4,i,\ell} <\infty$, $\sum_{\ell=1}^\infty \lambda_\ell
 \Delta_{i,\ell}^2 <\infty$ and $\sum_{\ell=1}^\infty \lambda_\ell
 \sigma_{4,i,\ell} <\infty$, for $i=2,\dots, k$. Then, $\sqrt n_i\left({\wbGa_i}-\bGa_1\right)\convdist \bU_i+ \tau_i^{1/2}\,\bDelta_i$, where  $\bDelta_i=\sum_{\ell\,\ge 1}\Delta_{i,\ell} \lambda_{\ell}
 \phi_\ell\otimes \phi_\ell$ and  $\bU_i$ a
  zero mean Gaussian random element with  
 covariance operator $\bUpsi_i$ given by (\ref{UPSIi1}) where $f_{im}$ are defined in (\ref{fcpc}).
   }

 \vskip0.2in

 \noi \textbf{Remark \ref{alt}.1.}  The conditions     $\sum_{\ell=1}^\infty \lambda_\ell |\Delta_{i,\ell}| <\infty$ and $\lambda_{\ell} (1+n^{-1/2}\Delta_{i,\ell})\ge 0$ ensure  that $\bDelta_i$ is a  self--adjoint trace operator and  that $\bGa_{i,n}$ is non--negative, respectively. Note that if, for all the populations,  the observations $X_{i,j}$ have a Gaussian distribution, then $f_{i \ell}\sim N(0,1)$, so   $\sigma^2_{4,i,\ell}=3$. This implies that the conditions   $\sum_{\ell=1}^\infty \lambda_\ell \Delta_{i,\ell}^2
 \sigma_{4,i,\ell} <\infty$, $\sum_{\ell=1}^\infty \lambda_\ell
 \Delta_{i,\ell}^2 <\infty$ and $\sum_{\ell=1}^\infty \lambda_\ell
 \sigma_{4,i,\ell} <\infty$ reduce to $\sum_{\ell=1}^\infty \lambda_\ell
 \Delta_{i,\ell}^2 <\infty$,  since $\sum_{\ell=1}^\infty \lambda_\ell<\infty$. Moreover, when considering root$-n$ local proportional alternatives, i.e., when $\Delta_{i,\ell}=\rho_i$, the  condition $\sum_{\ell=1}^\infty \lambda_\ell
 \Delta_{i,\ell}^2 <\infty$ is immediately fulfilled since $\bGa_1$ is a trace operator.
 
 \vskip0.2in
  Theorem  \ref{alt}.1 and Propositions \ref{alt}.1 and \ref{alt}.2, lead immediately to the asymptotic distribution of the test statistic $T_{k,n}$ under the  local alternatives studied above. We summarize this result   in Corollary \ref{alt}.1. 
 
 \vskip0.1in
 
 \noi  \textbf{Corollary \ref{alt}.1.} \textsl{Let $X_{i,1} ,\cdots,X_{i,n_i} $ for $i=1,\dots, k$ be independent observations
 from $k$ independent distributions  in $\itH$,  with mean $\mu_i$ and covariance
 operator $\bGa_i$ such that $\bGa_i=\bGa_{i,n}=\bGa_1+n^{-1/2}\bDelta_i$, for $i\ge 2$. Let $\wbGa_i$
  be the sample covariance operator of  the $i-$th
 population. Assume that the assumptions of Propositions \ref{alt}.1 or \ref{alt}.2 hold and denote $\bUpsi_{\weigh}=(\bUpsi_{\weigh, 1}, \dots, \bUpsi_{\weigh, k-1})$ where  $\bUpsi_{\weigh, i}$ is defined in (\ref{upsiweight}) with  $\bUpsi_{i}$  given in (\ref{UPSIi1}). Let  $\{\upsilon_\ell\}_{\ell\ge 1}$ be the orthonormal eigenfunctions of $\bUpsi_{\weigh}$ related to the eigenvalues $\{\theta_\ell\}_{\ell \ge 1}$ ordered in decreasing order and $\eta_\ell=\langle \bDelta^{(k-1)},v_\ell\rangle_{\itF^{k-1}}$. Then, we have that
   $$T_{k,n}=n \sum_{i=2}^k\|\wbGa_i -\wbGa_1 \|_\itF^2\convdist  \sum_{\ell\,\geq 1} \theta_\ell
   \left(Z_\ell+\frac{\eta_\ell}{\sqrt{\theta_\ell} }\right)^2\;,$$
   where $Z_\ell$ are independent, $Z_\ell \sim N(0,1)$. } 
 
 \vskip0.1in
  Under the local alternatives  $\bGa_{i,n}=\bGa_1+n^{-1/2}\bDelta_i$, for $i\ge 2$, and, in particular, under those given in Propositions \ref{alt}.1 and   \ref{alt}.2, similar arguments to those considered in the proof of Proposition 4 in Boente and Fraiman (2000) allow to show that, if $h=h_n\to 0$, the smooth estimator $\wbGa_{i,h}$ has the same asymptotic behaviour as  $\wbGa_i$, i.e., that $\sqrt{n_i}\| ( \wbGa_{i,h}- \bGa_{1,h})- (\wbGa_i-\bGa_1)\|_\itF\convprob 0$, where $\bGa_{1,h}$ is the smoothed covariance operator. On the other hand,  Proposition 3 in Boente and Fraiman (2000)  entails that $\sqrt{n}\|\bGa_{1,h}-\bGa_1)\|_\itF \to 0$ if, in addition, $n\,h \to 0$, the kernel $K$ has finite first moment and the covariance kernel $\gamma_1$ satisfies the following Lipschitz condition $|\gamma_1(t,u)-\gamma_1(t,t)|\le C |t-u|$, so that the asymptotic distribution of the statistic test $T_{k,n,h}=n \sum_{j=2}^k\|{\wbGa}_{j,h}-{\wbGa}_{1,h}\|_\itF^2$ is that given in Corollary \ref{alt}.1.

 \vskip0.1in
 \noi \textbf{Remark \ref{alt}.2.}
 Proportional alternatives of the form $\bGa_{i,n}=(1+\rho_i/\sqrt{n}) \bGa_1$   are obtained 
 taking $\Delta_{i,\ell}=\rho_i$ in Proposition \ref{alt}.2, so that  $\bDelta_i=\rho_i\bGa_1$. In this particular case,  we have that
 $$\left<\bDelta^{(k-1)},v_{i,i}\right>=\sum_{j=2}^{k}<\rho_j\bGa_1,\phi_i\otimes\phi_i>=\lambda_i \sum_{j=2}^{k}\rho_j$$
 and
 $$\left<\bDelta^{(k-1)},v_{i,j}\right>=\frac{1}{\sqrt 2}\sum_{j=2}^{k}<\rho_j\bGa_1,
 \phi_i\otimes\phi_j+\phi_j\otimes\phi_i>=0\,,$$
 where $\bGa_1=\sum \lambda_i \phi_i\otimes\phi_i$.  Moreover, if the processes are Gaussian, using Remark \ref{stest}.1, we get that the asymptotic distribution given in Theorem \ref{alt}.1, can be written as
 $W_k=2k^2\sum_{i\ge 1} \lambda_i^2(Z_i+ {\sum_{j=2}^{k}\rho_j}/({k\sqrt 2}))^2+ 2k^2\sum_{i\ge 1} \sum_{j\ge 1}\lambda_i\lambda_{i+j} Z_{i,j}^2$ and it   depends only on the eigenvalues of $\bGa_1$ different from zero.

\begin{sidewaysfigure}
 \small
 \begin{center}
 \vskip-0.3in
 \hskip-0.3in(a)\hskip3in (b)\\
 \vskip-0.3in
 \includegraphics[scale=0.4]{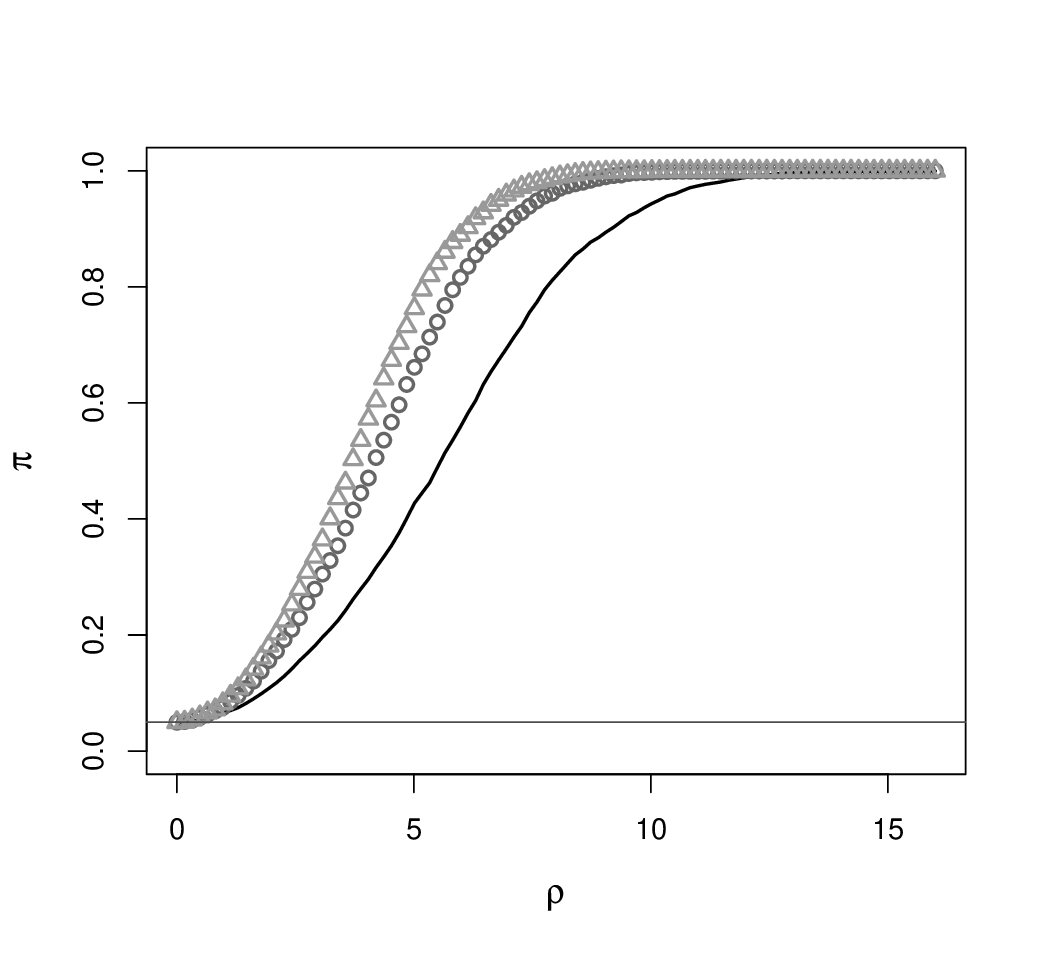} \hskip0.6in     \includegraphics[scale=0.4]{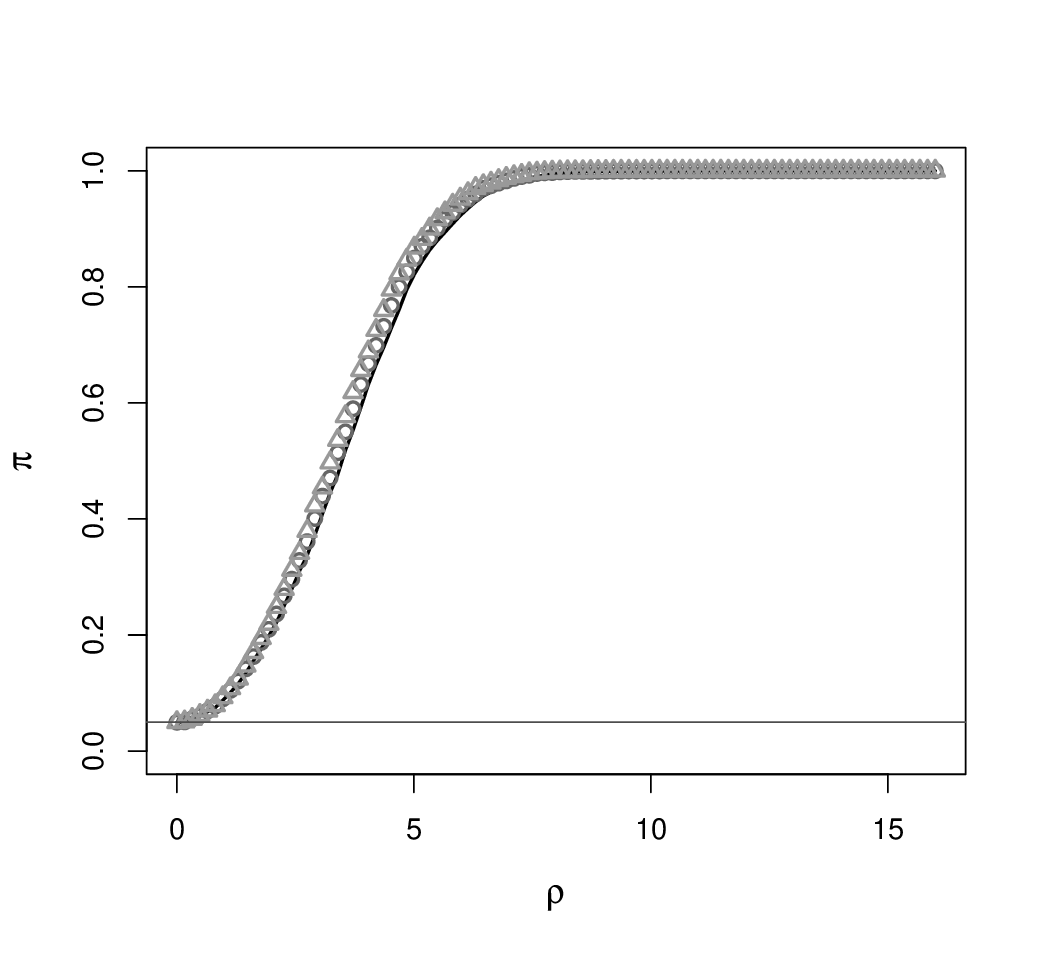}\\
 \vskip-0.1in
 \hskip-0.3in(c)\hskip3in (d)\\
 \vskip-0.3in
 $ $\hskip0.6in \includegraphics[scale=0.4]{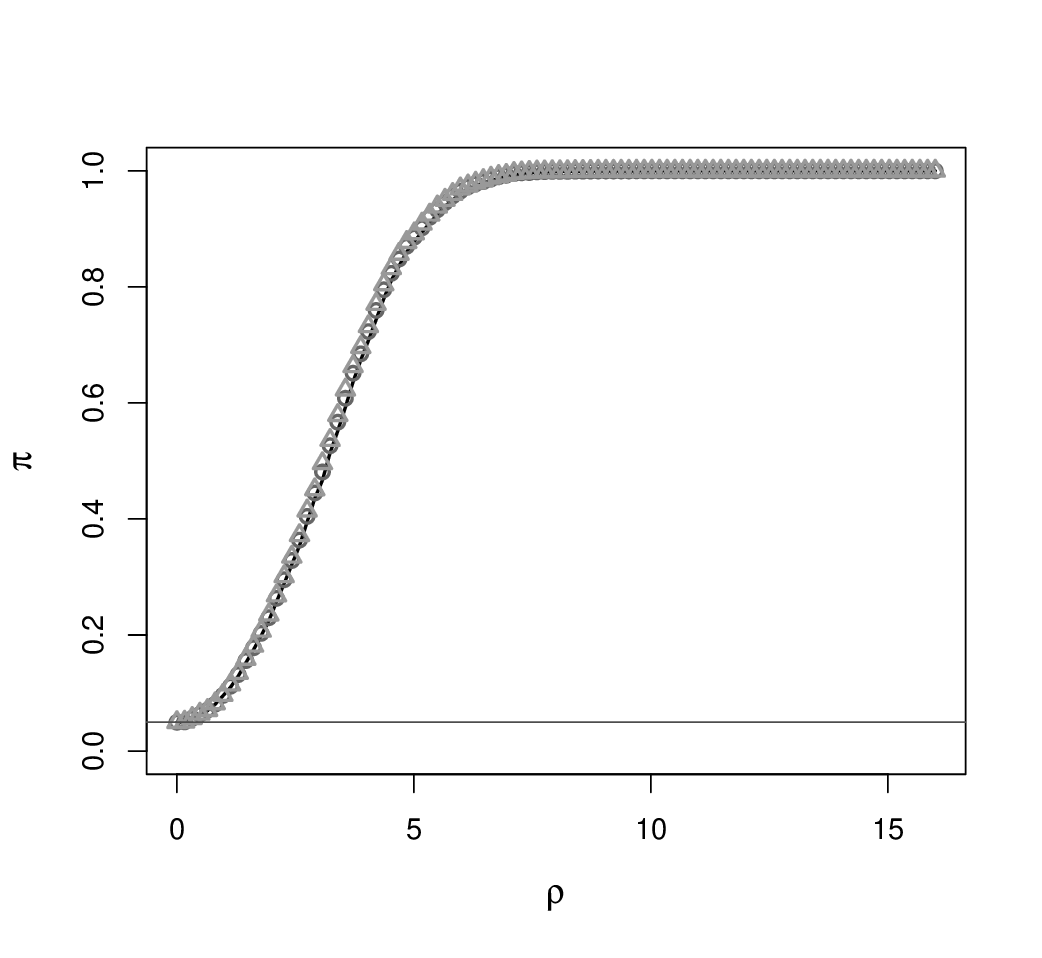}   \hskip-0.1in  \includegraphics[scale=0.56]{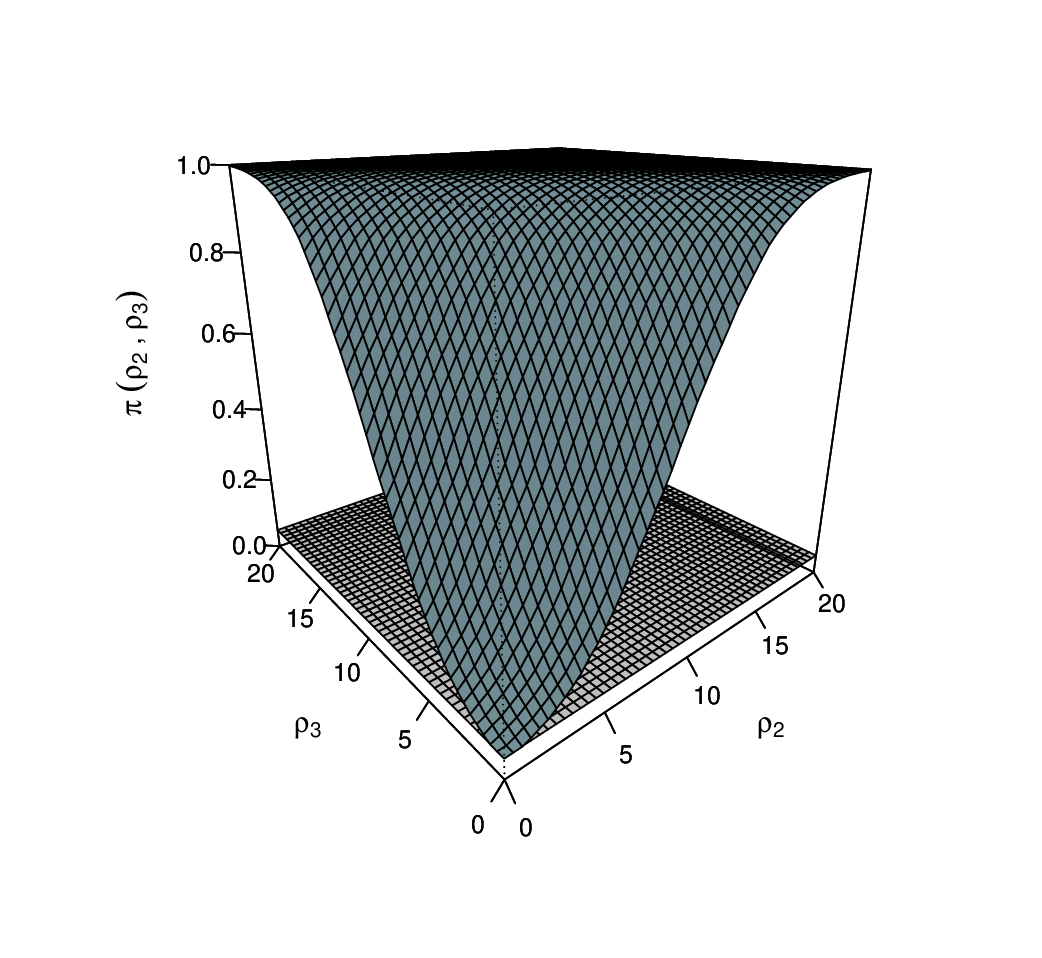}
  \end{center}
 \vskip-0.5in
 \caption{\label{fig:potencia}{\small
 Theoretical power for proportional Brownian motions. Figures a) to c) correspond to the situation  $\rho_1=\dots=\rho_k=\rho$, where 
 $k=2$ to $4$ in (a),   $k=5$ to $7$ in (b) and  $k=8$ to $10$ in (c).
 The solid lines correspond to $k=2,5,8$, the circles to $k=3,6,9$ and the triangles $k=4,7,10$. Figure (d) corresponds to $k=3$ and $\rho_2,\rho_3\in [0,20]$. }}
  \end{sidewaysfigure}

 Figure \ref{fig:potencia} contains the theoretical power computed using Monte Carlo, for
 different number of populations and alternatives, when the underlying processes are Brownian motions. In Figure  \ref{fig:potencia}a) to c), we choose  identical values of $\rho_j$, i.e., we considered the alternatives
 $ \Gamma_i=(1+{\rho}n^{-1/2})\Gamma_1 $, for $2\leq i\leq k$. On the other hand,   Figure  \ref{fig:potencia}d) corresponds to the three population
 situation and shows the surface plot of the theoretical power $\pi(\rho_2,\rho_3)$ when $ \Gamma_i=(1+ {\rho_i}n^{-1/2})\Gamma_1 $, for $2\leq i\leq 3$.
 
 To numerically compute the power, we have truncated $W_k$ as
 $$W_k=2k^2\sum_{i=1}^{20} \lambda_i^2(Z_i+ {\sum_{j=2}^{k}\rho_j}/({k\sqrt 2}))^2+ 2k^2\sum_{1\le i<j\le 20} \lambda_i\lambda_{j} Z_{i,j}^2\,.$$
 The value $20$ was chosen since the proportion of explained variance $\sum_{i=1}^{20}\lambda_i/\sum_{i\ge 1}\lambda_i$ is approximately $ 0.9898$.   Figure \ref{fig:potencia} (a) to (c)
 displays the theoretical  power $\pi(\rho)$ as a function of $\rho$ for  different values
 $\rho\in [0,10]$ and different number of populations. More precisely, Figure \ref{fig:potencia}  (a) corresponds to $k=2,3,4$,  (b) to $k=5,6,7$ and (c) to $k=8,9,10$.  The solid lines correspond to $k=2,5,8$, the circles to $k=3,6,9$ and the triangles $k=4,7,10$. On the other hand,    Figure \ref{fig:potencia}(d) provides a surface plot for the theoretical power  $\pi(\rho_2,\rho_3)$ when $\rho_i\in [0,20]$ for $i=2,3$. The horizontal gray line in a) to c) and the horizontal gray plane in d) correspond to the level 0.05. These plots show that the test improves its performance considerably when $k=3$ populations are compared instead of two populations. Besides, the power is quite stable for values of $k$ larger than $5$ and for the proportional alternatives considered it shows an important detection capability,  when $k\ge 4$.

 \vskip0.1in
  When the stronger condition $\sup_{n\ge 1} \esp \|X_i^{(n)}\|^{4+\delta}<\infty$ holds, Theorem 2.1  in Gaines \textsl{et al.} (2011) 
    together with Theorem \ref{alt}.1 lead immediately to the asymptotic distribution of  test statistic $T_{k,n}$ under root$-n$ local alternatives as stated  in Proposition \ref{alt}.3. 
 
 \vskip0.1in
 
 \noi  \textbf{Proposition \ref{alt}.3.} \textsl{Let $X_{i,1}^{(n)} ,\cdots,X_{i,n_i}^{(n)} $ for $i=1,\dots, k$ be independent observations
 from $k$ independent distributions  in $\itH$,  with  covariance operators $\bGa_i$ such that, for $i\ge 1$, $\bGa_i=\bGa_{i,n}=\bGa_1+n^{-1/2}\bDelta_i$, where  $\bDelta_i$ is a trace operator and $\sup_{n\ge 1} \esp \|X_i^{(n)}\|^{4+\delta}<\infty$, with $X_{i,j}^{(n)}\sim X_{i}^{(n)}$. Assume that, for $i\ge 2$,   the covariance operator of $\bY_{i}^{(n)}=(X_{i}^{(n)}-\esp(X_{i}^{(n)}))\otimes (X_{i}^{(n)}-\esp(X_{i}^{(n)}))$  converges to an operator $\bUpsi_i$ in trace norm. Denote $\bUpsi_{\weigh}=(\bUpsi_{\weigh, 1}, \dots, \bUpsi_{\weigh, k-1})$ where  $\bUpsi_{\weigh, i}$ is defined in (\ref{upsiweight}) and $\bUpsi_{1}$ is given in (\ref{UPSIi1}). Let  $\{\upsilon_\ell\}_{\ell\ge 1}$ be the orthonormal eigenfunctions of $\bUpsi_{\weigh}$ related to the eigenvalues $\{\theta_\ell\}_{\ell \ge 1}$ ordered in decreasing order and $\eta_\ell=\langle \bDelta^{(k-1)},v_\ell\rangle_{\itF^{k-1}}$. Then if  $\wbGa_i$
  stands for the sample covariance operator of  the $i-$th
 population, we have that
   $$T_{k,n}=n \sum_{i=2}^k\|\wbGa_i -\wbGa_1 \|_\itF^2\convdist  \sum_{\ell\,\geq 1} \theta_\ell
   \left(Z_\ell+\frac{\eta_\ell}{\sqrt{\theta_\ell} }\right)^2\;,$$
   where $Z_\ell$ are independent, $Z_\ell \sim N(0,1)$. } 
 
 \vskip0.1in
  It is worth noting that if $X_{i,j} $ satisfy (\ref{contaminado}) and $ \esp \|W_{i}\|^{4+\delta}<\infty$ and  $\esp \|R_{i}\|^{4+\delta}<\infty$ the proof of Proposition \ref{alt}.1 is a consequence of  Theorem 2.1  in Gaines \textsl{et al.} (2011).  Similarly, if  $\Delta_{i,\ell}\ge 0$ and $ \esp \|X_i^{(1)}\|^{4+\delta}<\infty$ the proof of Proposition \ref{alt}.1  can also be derived from Theorem 2.1  in Gaines \textsl{et al.} (2011)  through straightforward calculations. However, in both cases, we prefer to avoid imposing higher moment conditions and/or to consider more general alternatives and for that reason we have included their proof   in the Appendix.

 \section{Bootstrap calibration}{\label{bootstrap}}
 The asymptotic null behaviour derived in Section \ref{stest} motivate the use of the bootstrap methods, due the fact that  the asymptotic distribution obtained in (\ref{distkpob})  depends on the   unknown eigenvalues $\theta_\ell $.
 For that reason, we will consider a general bootstrap method to approximate the distribution of the test which can be described as follows.
 \begin{itemize}
 \item[] \textbf{Step  1.} For $1\le i\le k$, and given the sample $X_{i,1},\cdots,X_{i,n_i}$, let  $\wbUps_{i}$ be consistent estimators of $\bUpsi_i$. Define
 $\wbUps_{\weigh}=(\wbUps_{\weigh, 1},\dots, \wbUps_{\weigh, k-1})$ where
 $$
 \wbUps_{\weigh, i}(u_1,\dots,u_{k-1})= \frac{1}{\wtau_{i+1}}\wbUps_{i+1}(u_1)+\frac{1}{\wtau_1}\wbUps_1 \left(\sum_{i=1}^{k-1}u_i\right)\,,
 $$
 and $\wtau_i=n_i/\sum_{s=1}^k n_s$.
 In particular, if $k=2$, $\wbUps_{\weigh}=
 \wtau_1^{\,-1}\wbUps_1+\wtau_2^{\,-1}\wbUps_2$ with $\wtau_i=n_i/({n_1+n_2})$.
 \item[]\textbf{Step  2.} For $1\le \ell \le q_n$ denote by  $\wtheta_\ell $  the positive eigenvalues of  $\wbUps_{\weigh}$.
 \item[]\textbf{Step  3.} Generate $Z^*_1,\dots,Z^*_{q_n}$ i.i.d. such that $Z^*_i \sim N(0,1)$ and let ${\itU}^*_{n}=\sum_{j=1}^{q_n}\wtheta_j {Z^*_{j}}^2$.
 \item[]\textbf{Step  4.}  Repeat \textbf{Step  3}  $N_{\boot}$ times, to get $N_{\boot}$ values of ${\itU}_{nr}^*$ for $1\leq r\leq N_{\boot}$.
 \end{itemize}
 The $(1-\alpha)-$quantile of the asymptotic distribution of $T_{k,n}$ can be approximated by the $(1-\alpha)-$quantile of
 the empirical distribution of ${\itU}_{nr}^*$ for $1\leq r\leq N_{\boot}$.  The $p-$value can be estimated by $\widehat{p}=s/{N_{\boot}}$ where $s$ is the number of ${\itU}^*_{nr}$ which are larger or equal than the observed value of $T_{k,n}$.
 
 \vskip0.1in
 \noi \textbf{Remark \ref{bootstrap}.1.}
  Note that this procedure depends only on the asymptotic distribution of $\wbGa_i$. For the sample covariance estimator, the covariance operator $\bUpsi_i$ to be estimated in \textbf{Step 1} is given in (\ref{varasintgamai}). Assume  that all the populations have a Gaussian distribution, then    $\bUpsi_i$ can be  estimated using the eigenvalues and eigenfunctions of the sample covariance, since  $f_{ij}$ are independent and $f_{ij}\sim N(0,1)$. For non Gaussian samples, $\bUpsi_i$ can be  estimated noticing that it is the covariance operator of $\bY_i=(X_{i}-\mu_i)\otimes (X_{i}-\mu_i)$.
 When considering other asymptotically normally estimators of $\bGa_i$, such as  the smoothed estimators $\wbGa_{i,h}$ for $L^2({\itI})$ trajectories, the estimators need to be adapted.

  \vskip0.2in
 Note that the space of covariance operators of random elements on $\itH$ is a Hilbert space with the inner product defined in $\itF$. Hence, the covariance of any estimate of the covariance operator  is also an element of a Hilbert space, which we denote as $\itG$. Then,  for instance, $\bUpsi_i$ and $\wbUps_{i}$ in \textbf{Step 1} belong to $\itG$, while $\wbUps_{\weigh}$ and $\bUpsi_{\weigh}$ are  random elements of the product Hilbert space $\itG^{k-1}$ with norm denoted as $\|\cdot\|_{\itG^{k-1}}$.
 
 \vskip0.1in
 The following theorem entails the validity of the bootstrap calibration method.
 It states that  the bootstrap distribution of ${\itU}^*_n$
 converges to the asymptotic null distribution of $T_n$.   This fact
 ensures that the asymptotic significance level of the test based on
 the bootstrap critical value is indeed $\alpha$   and that the bootstrap test leads to a consistent test.  
 
 \vskip0.1in
 \noi \textbf{Theorem \ref{bootstrap}.1.} \textsl{
 Let $q_n$ such that $q_n/\sqrt n\to 0$ and  $\wtX_n=(X_{1,1},\cdots,X_{1,n_1},\dots, X_{k,1},\cdots,X_{k,n_k})$. Denote by $F_{{\itU}^*_n\vert \wtX_n}(\cdot)=\prob({\itU}^*_n\leq \cdot\;\vert \wtX_n)$ and by $F_{\itU}$ the distribution function of ${\itU}=\sum_{\ell\,\geq 1} \theta_\ell  Z^2_\ell$, where $Z_{\ell}$ are i.i.d. and $Z_{\ell}\sim N(0,1)$. Assume that  $\esp(\|X_{i}\|^4)<\infty$  
 and ${n_i}/n\to \tau_i$ with $\tau_i\in (0,1)$ and $n=\sum_{i=1}^k n_i$.
  Then,  if $\sqrt n \|\wbUps_{\weigh}-\bUpsi_{\weigh}\|_{\itG^{k-1}}=O_{\prob}(1)$,  we have that
 $\rho_{\mbox{\footnotesize\sc k}}( F_{{\itU}^*_n\vert \wtX_n}, F_{\itU})\convprob 0\;,$
 where   $\rho_{\mbox{\footnotesize\sc k}}(F,G)$ stands for the Kolmogorov distance between the distribution functions $F$ and $G$.
 }

 \section{Simulation study}\label{simu}
  This section contains the results of  two simulation studies carried on
 with $k=2$ and $k=3$ populations and designed to illustrate the finite--sample performance of the bootstrap test procedure described in Section  \ref{bootstrap},
 under the null hypothesis and under different alternatives. 
 In all cases, we generate $NR=1000$ samples of size $n_i$, $1\le i\le k$ and each trajectory was observed at  $m=100$ equidistant points in the interval $[0,1]$. To analyse the dependence on the sample size, we choose $n_i =50,  100$ and $ 200$, for $1\le i\le k$ which allows to study the behaviour of the test in terms of level approximation as well as power performance depending on the sample size. In all tables, we report the observed frequency of rejections over replications with nominal level $\alpha=0.05$.
 
 \subsection{Simulation settings}
  
 Under the null hypothesis, we consider infinite--dimensional processes generating   independent centred  Brownian motion processes,  denoted from now on as ${\cal{BW}}(0,1)$. On the other hand, to check the test power performance, we consider  root$-n$ local alternatives. 
 To be more precise, when comparing two populations, 
 we generate independent observations  $X_{1,j}\sim X_1$, $1\le j\le n_1$, and $X_{2,j}\sim X_2$, $1\le j\le n_2$, such that $X_{1}\sim {\cal{BW}}(0,1)$ and $X_{2}\sim W_1+\delta_n\, W_2^2$, where $W_1$ and $W_2$ are independent $W_i\sim {\cal{BW}}(0,1)$, $i=1,2$ and $\delta_n={\rho}n^{-1/4}$ with $n= n_1+n_2$. The situation $\rho=0$ corresponds to the null hypothesis, while to study the test power  we choose  $\delta_n={\rho}n^{-1/4}$ with $n= n_1+n_2$ and 
 $\rho$ taking   values from 1 to 10.  
 Note that for set of alternatives, the covariance operator of $X_{2,1}$ equals $\bGa_2= \bGa_1 + \rho^2\,n^{-1/2}\, \bDelta$, where $\bDelta$ is the covariance operator of $W_2^2$. These alternatives   correspond to the local alternatives studied in Proposition \ref{alt}.1.
  
 On the other hand, for  the  three populations case, we consider a proportional model taking independent observations  $X_{i,j}\sim X_i$, $1\le j\le n_i$, $1\le i\le k$, such that $X_{1}\sim {\cal{BW}}(0,1)$, while
 $X_{i}\sim (1+\delta_n)^{1/2}{\cal{BW}}(0,1)$, for $i=2,3$, where $\delta_n={\rho}n^{-1/2}$ with $n=\sum_{i=1}^3 n_i$. The parameter $\rho$ takes  values   on an equidistant grid of points between 0 and 20 of size 11. In this case,  the covariance operators of $X_{2 }$ and $X_{3 }$ equal $\bGa_2=\bGa_3= (1+{\rho}n^{-1/2})\bGa_1 $ corresponding to the proportional alternatives described in Remark \ref{alt}.2.
 
 \subsection{The  testing procedures}
 We study the behaviour of  the test based on $T_{k,n}$ defined in (\ref{tkn}) using the bootstrap calibration described in Section \ref{bootstrap} with $N_{\boot} =5000$ bootstrap replications.  To perform the bootstrap calibration,   we  project the centred  data onto the   $M$  largest principal components of the pooled sample covariance matrix $n^{-1}\sum n_i \wbGa_i$. 
 We then estimate the covariance operator  $\wbUps_{\weigh}$ through a finite dimensional matrix. To evaluate the dependence on the number of principal components chosen, we select different values of $M$ as $M=  3, 10,  20, 30 $. Note that, in this situation, the value $q_n$ used in \textbf{Step 2} equals $q_n=M(M+1)/2$. The percentage of total variance explained  by the selected number of principal components is reported in Table  \ref{tab:percent}, while  the frequencies of rejection  at the 5$\%$ level, for $k=2$ and $k=3$, are given in Tables  \ref{tab:power_k2} and \ref{tab:power_k3}, respectively.  
 
 \begin{table}[ht!]
     \centering
    \small
 \begin{tabular}{c|c|cccc|cccc|cccc|}
  \hline
  & & \multicolumn{4}{c|}{$n_i=50$} & \multicolumn{4}{c|}{$n_i=100$} & \multicolumn{4}{c|}{$n_i=200$}\\\hline
  $k$ & $\rho$& \multicolumn{4}{c|}{$M$}   & \multicolumn{4}{c|}{$M$}  & \multicolumn{4}{c|}{$M$} \\\hline
  & & 3 & 10 & 20 & 30 & 3 & 10 & 20 & 30 & 3 & 10 & 20 & 30 \\
   \hline
 $2$ & 0 & 0.935 & 0.982 & 0.993 & 0.996 & 0.934 & 0.981 & 0.991 & 0.995 & 0.934 & 0.980 & 0.991 & 0.994 \\
 $3$ & 0 & 0.934 & 0.981 & 0.991 & 0.995 & 0.934 & 0.980 & 0.991 & 0.995 & 0.934 & 0.981 & 0.991 & 0.994 \\ 
    \hline
 \end{tabular}
 \caption{\small \label{tab:percent} Percentage of the total variance explained by the first $M$ principal components.}
 \end{table}
 
 Taking into account the fact that, under the null hypothesis, the processes are Gaussian,  
     Remark \ref{stest}.1.d) entails that  $\theta_{i,i}=2k^2\lambda_i^2$ and $\theta_{i,j}=2k^2\lambda_i\lambda_j$, for $i<j$. Then, from the eigenvalues $\wlam_{\ell}$ of  the pooled sample covariance matrix,  one may easily  provide estimators $\wtheta_j$ of $\theta_j$ to replace those considered in \textbf{Step 2}. This approximation is denoted as \textsl{Gaussian} in Tables \ref{tab:power_k2} and \ref{tab:power_k3} and was computed using the fact 
 that the trajectories were generated over a  grid of $100$ points for all the sample sizes leading to at most $100$ non--null values $\wlam_{\ell}$.

 We also compare the behaviour of our test statistic with the permutation test introduced in Pigoli \textsl{et al.} (2014) when $k=2$. Our choice for the permutation test is based on the numerical study reported in Pigoli \textsl{et al.} (2014), where it is shown that the permutation test provides a good competitor to the tests introduced in Panaretos \textsl{et al.} (2010) and Fremdt \textsl{et al.} (2013). We perform the permutation test taking the same   discrepancy measure between covariance operators   used for $T_{k,n}$, i.e., $d(\bGa_1, \bGa_2)=\|\bGa_1 -\bGa_2 \|_\itF$. The obtained results when using $N_{\perm}=1000$ and $5000$ random permutations are given in  Table \ref{tab:power_k2_perm}. In the case of $k=3$ populations, a permutation  test was also considered taking as test statistic  $D=d(\wbGa_1, \wbGa_2,\wbGa_3)=\|\wbGa_2 -\wbGa_1 \|_\itF^2+\|\wbGa_3 -\wbGa_1 \|_\itF^2+\|\wbGa_3 -\wbGa_2 \|_\itF^2$. As in Pigoli \textsl{et al.} (2014), we first center the samples   using the sample mean and then, we consider $N_{\perm}$ random permutations of the labels ${1, 2,3}$ on the centred sample curves. For each permutation $j$, we compute $D_j=d(\wbGa_1^{(j)}, \wbGa_2^{(j)},\wbGa_3^{(j)})$, for $j=1,\dots, N_{\perm}$, where $\wbGa_i^{(j)}$ is  the sample covariance operator of the group indexed with label $i$ in the given permutation. As in the two population case, the $p-$value of the test is the proportion of $D_j$  which are greater than or equal than $D$. Table \ref{tab:power_k3_perm_mod} reports the obtained frequencies of rejection. We also used this approach taking as test statistic $D^{\star}=d(\wbGa_1, \wbGa_2,\wbGa_3)=\|\wbGa_2 -\wbGa_1 \|_\itF^2+\|\wbGa_3 -\wbGa_1 \|_\itF^2$, which corresponds to $T_{k,n}$, but is not invariant by permutation of the labels. The results for $D^{\star}$ are similar to those obtained for $D$ and are not reported here.
 
 In all tables, we denote as  $\phi_{\boot,M}$, for $M=3,10,20$ and $30$ the bootstrap calibration of $T_{k,n}$ computed using $M$ principal components, $\phi_{\boot,\gaus}$ the bootstrap calibration of $T_{k,n}$ computed using the Gaussian approximation for $\theta_{i,j}$ and $\phi_{\perm,N_{\perm}}$, for $N_{\perm}=1000$ and $5000$, the permutation test computed using $N_{\perm}$ random permutations.
 
 \subsection{Simulation results}
 With respect to the bootstrap calibration described in Section \ref{bootstrap} for the test based on $T_{k,n}$,  Tables \ref{tab:power_k2} and \ref{tab:power_k3} show the
 improvement attained in level when the Gaussian  approximation  is used, both for   $k=2$ and $k=3$  populations.  Also, when we project the data on the first $M$ principal components,  the empirical size of the   test based on the bootstrap  calibration is quite close to 
 the nominal one. To analyse the significance of the empirical size,  we study if the empirical size  is significantly 
 different from the nominal level $\alpha=0.05$. To be more precise, for a test $\phi_n$ based on a statistic $T_n$, let $\pi$ be such that $\pi_{H_0}(\phi_n)\convprob \pi$. Then, using the central limit theorem,  
 the hypothesis  $H_{0,\pi}: \pi=\alpha$ is rejected at level $\gamma$ versus  $H_{1,\pi}:\pi\ne\alpha$  
 if $\pi_{H_0}(\phi_n)\notin [a_1(\alpha),a_2(\alpha)]$ where $a_j(\alpha)=\alpha+ (-1)^j z_{\gamma/2}\,\{\alpha(1-\alpha)/NR\}^{1/2}$, 
 $j=1,2$. If $H_{0,\pi}: \pi=\alpha=0.05$  is not rejected, the testing procedure based on $T_n$ is 
 considered accurate, while if $\pi_{H_0}(\phi_n)<a_1(\alpha)$  the testing procedure is  conservative and when $\pi_{H_0}(\phi_n)>a_2(\alpha)$ the test is liberal. In all the considered situations for $k=2$,  the test is accurate with significance level $\gamma=0.01$. On the other hand, for $k=3$ populations, the test is liberal only when $n_1=n_2=n_3=50$ and $M=3$ or $10$, in all other situations the test is accurate, so in almost all considered situations the proposed method has a  good level performance. 
  
 Regarding the behaviour under the alternative,  the  bootstrap test  detects  
 the considered alternatives  for different values of $M$ and also when using the Gaussian approximation to the eigenvalues $\theta_{i,j}$. As expected, the observed frequencies of rejection converge to one as $\rho$ increases. Since local
 alternatives are taken, the power is almost similar for all choices of sample sizes and shows the tests capability to detect
 the selected local alternatives.   However, it is worth noting that the test shows a slower power convergence  for $k=2$ and $n_1=n_2=50$.

 \begin{sidewaystable} 
 \small    \centering
 \begin{tabular}{c|cccc|c||cccc|c||cccc|c|}
  \hline
  & \multicolumn{5}{c||}{$n_1=n_2=50$} & \multicolumn{5}{c||}{$n_1=n_2=100$} & \multicolumn{5}{c|}{$n_1=n_2=200$}\\\hline
  $\rho$ &  $\phi_{\boot,3}$ &  $\phi_{\boot,10}$ &  $\phi_{\boot,20}$ &  $\phi_{\boot,30}$ & $\phi_{\boot,\gaus}$ & $\phi_{\boot,3}$ &  $\phi_{\boot,10}$ &  $\phi_{\boot,20}$ &  $\phi_{\boot,30}$ & $\phi_{\boot,\gaus}$ & $\phi_{\boot,3}$ &  $\phi_{\boot,10}$ &  $\phi_{\boot,20}$ &  $\phi_{\boot,30}$ & $\phi_{\boot,\gaus}$ \\
   \hline
   0 & 0.068 & 0.064 & 0.063 & 0.061 & 0.048 & 0.066 & 0.065 & 0.061 & 0.060 & 0.050 & 0.054 & 0.052 & 0.052 & 0.051 & 0.040 \\ 
   1 & 0.074 & 0.070 & 0.068 & 0.067 & 0.066 & 0.083 & 0.082 & 0.081 & 0.079 & 0.064 & 0.081 & 0.078 & 0.076 & 0.075 & 0.059 \\ 
   2 & 0.234 & 0.218 & 0.215 & 0.208 & 0.305 & 0.315 & 0.299 & 0.296 & 0.290 & 0.333 & 0.356 & 0.348 & 0.343 & 0.337 & 0.355 \\ 
   3 & 0.536 & 0.512 & 0.498 & 0.490 & 0.721 & 0.694 & 0.681 & 0.671 & 0.666 & 0.801 & 0.851 & 0.845 & 0.840 & 0.837 & 0.891 \\ 
   4 & 0.722 & 0.699 & 0.689 & 0.682 & 0.911 & 0.895 & 0.890 & 0.885 & 0.882 & 0.975 & 0.992 & 0.990 & 0.989 & 0.988 & 0.999 \\ 
   5 & 0.814 & 0.796 & 0.788 & 0.785 & 0.979 & 0.948 & 0.942 & 0.941 & 0.940 & 0.998 & 0.999 & 0.999 & 0.999 & 0.998 & 1.000 \\ 
   6 & 0.851 & 0.836 & 0.829 & 0.821 & 0.997 & 0.959 & 0.957 & 0.956 & 0.956 & 1.000 & 0.999 & 0.999 & 0.999 & 0.999 & 1.000 \\ 
   7 & 0.864 & 0.853 & 0.847 & 0.840 & 0.999 & 0.967 & 0.962 & 0.961 & 0.961 & 1.000 & 0.999 & 0.999 & 0.999 & 0.999 & 1.000 \\ 
   8 & 0.872 & 0.857 & 0.853 & 0.850 & 1.000 & 0.972 & 0.969 & 0.965 & 0.964 & 1.000 & 0.999 & 0.999 & 0.999 & 0.999 & 1.000 \\ 
   10 & 0.873 & 0.864 & 0.859 & 0.857 & 1.000 & 0.973 & 0.971 & 0.968 & 0.967 & 1.000 & 0.999 & 0.999 & 0.999 & 0.999 & 1.000 \\ 
   
 \hline
 \end{tabular}
 \caption{\small \label{tab:power_k2} Frequency of rejection for the bootstrap test  $\phi_{\boot,M}$ when   $M=3,10,20$ and $30$ principal components are used for different sample sizes and two populations are compared. The  column labelled $\phi_{\boot,\gaus}$ reports the frequencies obtained when the eigenvalues $\theta_\ell$ are estimated using that the processes are Gaussian as described in Remark \ref{stest}.1.d). The alternatives considered are $X_{1 }\sim {\cal{BW}}(0,1)$ while
 $X_{2 }\sim W_1+\delta n^{-1/4} W_2^2$, where $W_j\sim {\cal{BW}}(0,1)$ are independent of each other and $\delta_n={\rho}n^{-1/4}$ with $n= n_1+n_2$. }
 
 \vskip0.1in

 \begin{tabular}{c|cccc|c||cccc|c||cccc|c|}
  \hline
  & \multicolumn{5}{c||}{$n_1=n_2=n_3=50$} & \multicolumn{5}{c||}{$n_1=n_2=n_3=100$} & \multicolumn{5}{c|}{$n_1=n_2=n_3=200$}\\\hline
  $\rho$ &  $\phi_{\boot,3}$ &  $\phi_{\boot,10}$ &  $\phi_{\boot,20}$ &  $\phi_{\boot,30}$ & $\phi_{\boot,\gaus}$ & $\phi_{\boot,3}$ &  $\phi_{\boot,10}$ &  $\phi_{\boot,20}$ &  $\phi_{\boot,30}$ & $\phi_{\boot,\gaus}$ & $\phi_{\boot,3}$ &  $\phi_{\boot,10}$ &  $\phi_{\boot,20}$ &  $\phi_{\boot,30}$ & $\phi_{\boot,\gaus}$ \\
   \hline
   0 & 0.071 & 0.069 & 0.063 & 0.062 & 0.065 & 0.066 & 0.062 & 0.059 & 0.058 & 0.066 & 0.064 & 0.058 & 0.057 & 0.054 & 0.058 \\ 
   2 & 0.126 & 0.118 & 0.118 & 0.115 & 0.099 & 0.139 & 0.135 & 0.134 & 0.134 & 0.108 & 0.114 & 0.108 & 0.107 & 0.104 & 0.094 \\ 
   4 & 0.272 & 0.258 & 0.255 & 0.251 & 0.216 & 0.285 & 0.276 & 0.275 & 0.272 & 0.243 & 0.294 & 0.280 & 0.273 & 0.271 & 0.242 \\ 
   6 & 0.437 & 0.426 & 0.422 & 0.418 & 0.369 & 0.476 & 0.460 & 0.455 & 0.453 & 0.413 & 0.496 & 0.481 & 0.476 & 0.472 & 0.454 \\ 
   8 & 0.623 & 0.610 & 0.604 & 0.602 & 0.541 & 0.663 & 0.656 & 0.652 & 0.646 & 0.606 & 0.695 & 0.683 & 0.680 & 0.678 & 0.668 \\ 
   10 & 0.746 & 0.733 & 0.728 & 0.725 & 0.686 & 0.798 & 0.793 & 0.792 & 0.790 & 0.760 & 0.843 & 0.832 & 0.829 & 0.828 & 0.818 \\ 
   12 & 0.843 & 0.834 & 0.830 & 0.829 & 0.791 & 0.899 & 0.896 & 0.890 & 0.889 & 0.875 & 0.926 & 0.922 & 0.917 & 0.914 & 0.908 \\ 
   14 & 0.909 & 0.904 & 0.901 & 0.899 & 0.867 & 0.951 & 0.949 & 0.944 & 0.944 & 0.943 & 0.976 & 0.972 & 0.967 & 0.966 & 0.968 \\ 
   16 & 0.945 & 0.940 & 0.938 & 0.937 & 0.928 & 0.973 & 0.971 & 0.970 & 0.969 & 0.969 & 0.993 & 0.993 & 0.992 & 0.991 & 0.990 \\ 
   18 & 0.971 & 0.968 & 0.968 & 0.967 & 0.961 & 0.991 & 0.989 & 0.989 & 0.988 & 0.985 & 0.999 & 0.999 & 0.999 & 0.999 & 0.999 \\ 
   20 & 0.986 & 0.985 & 0.985 & 0.985 & 0.982 & 0.995 & 0.995 & 0.995 & 0.995 & 0.994 & 1.000 & 1.000 & 1.000 & 1.000 & 1.000 \\ 
    \hline
 \end{tabular}
 \caption{\small \label{tab:power_k3} Frequency of rejection for the bootstrap test  $\phi_{\boot,M}$  when   $M=3,10,20$ and $30$ principal components are used for different sample sizes and $k=3$ populations are considered, when   $X_{1}\sim {\cal{BW}}(0,1)$ and $X_{i}\sim (1+\delta_n)^{1/2}{\cal{BW}}(0,1)$, for $i=2,3$, where $\delta_n={\rho}n^{-1/2}$ with $n=\sum_{i=1}^3 n_i$.  The column labelled $\phi_{\boot,\gaus}$ reports the frequencies obtained when the eigenvalues $\theta_\ell$ are estimated using that the processes are Gaussian as described in Remark \ref{stest}.1.d).}
 
 \end{sidewaystable}
 
 \normalsize
  
\begin{sidewaysfigure}
 \centering
 $n_1=n_2=50$ \hskip1.7in $n_1=n_2=100$ \hskip1.7in $n_1=n_2=200$\\ 
 \includegraphics[scale=0.4]{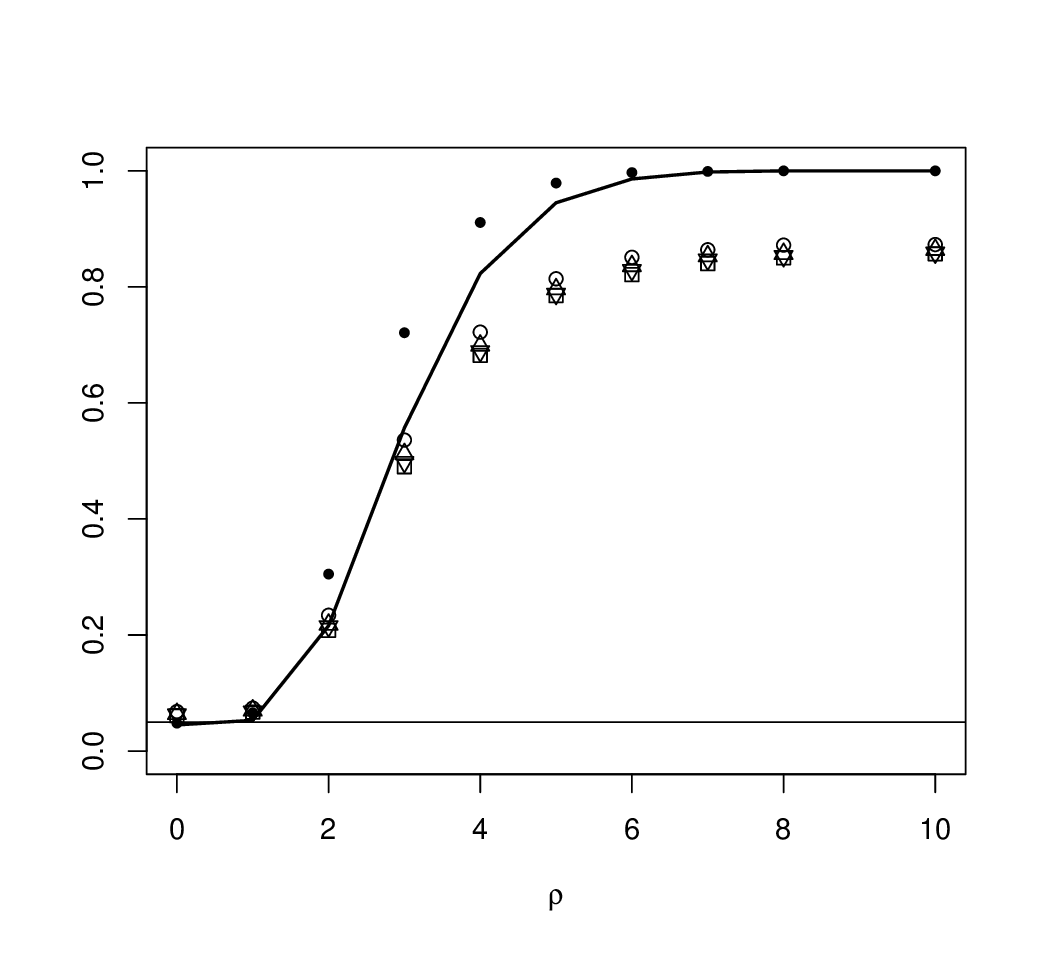} \hspace{-0.5cm}  \includegraphics[scale=0.4]{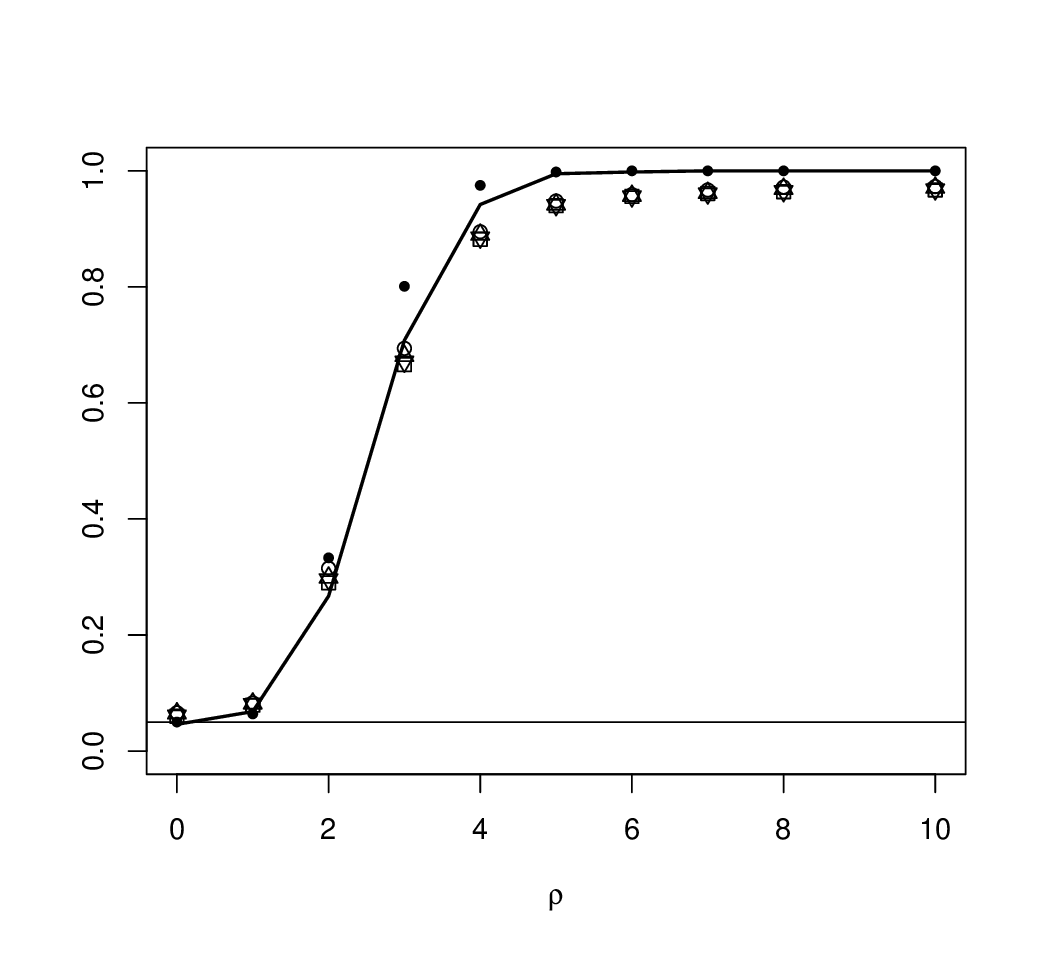}  \hspace{-0.5cm}\includegraphics[scale=0.4]{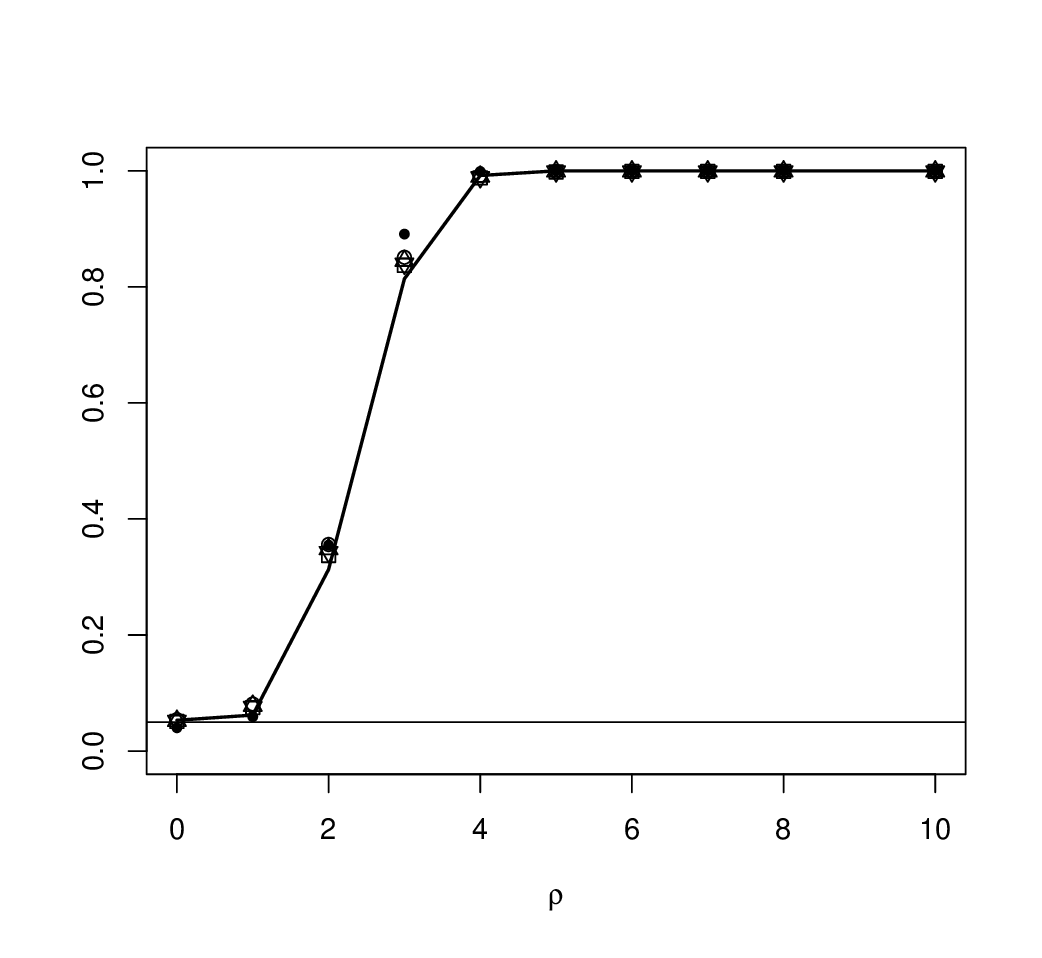} 
  \vspace{-0.5cm}
  \caption{\small \label{fig:power_k2}  Frequency of rejection when $k=2$ for the bootstrap test $\phi_{\boot,M}$, $\phi_{\boot,\gaus}$, and  $\phi_{\perm,5000}$. The solid line corresponds to $\phi_{\perm,5000}$, the filled circles to $\phi_{\boot,\gaus}$, while the circles, upper, lower triangles and the square correspond to $\phi_{\boot,M}$, with $M=3,10,20$ and $30$, respectively. }
 
 \vskip0.1in
 $n_1=n_2=n_3=50$ \hskip1.4in $n_1=n_2=n_3=100$ \hskip1.4in $n_1=n_2=n_3=200$\\ 
 \includegraphics[scale=0.4]{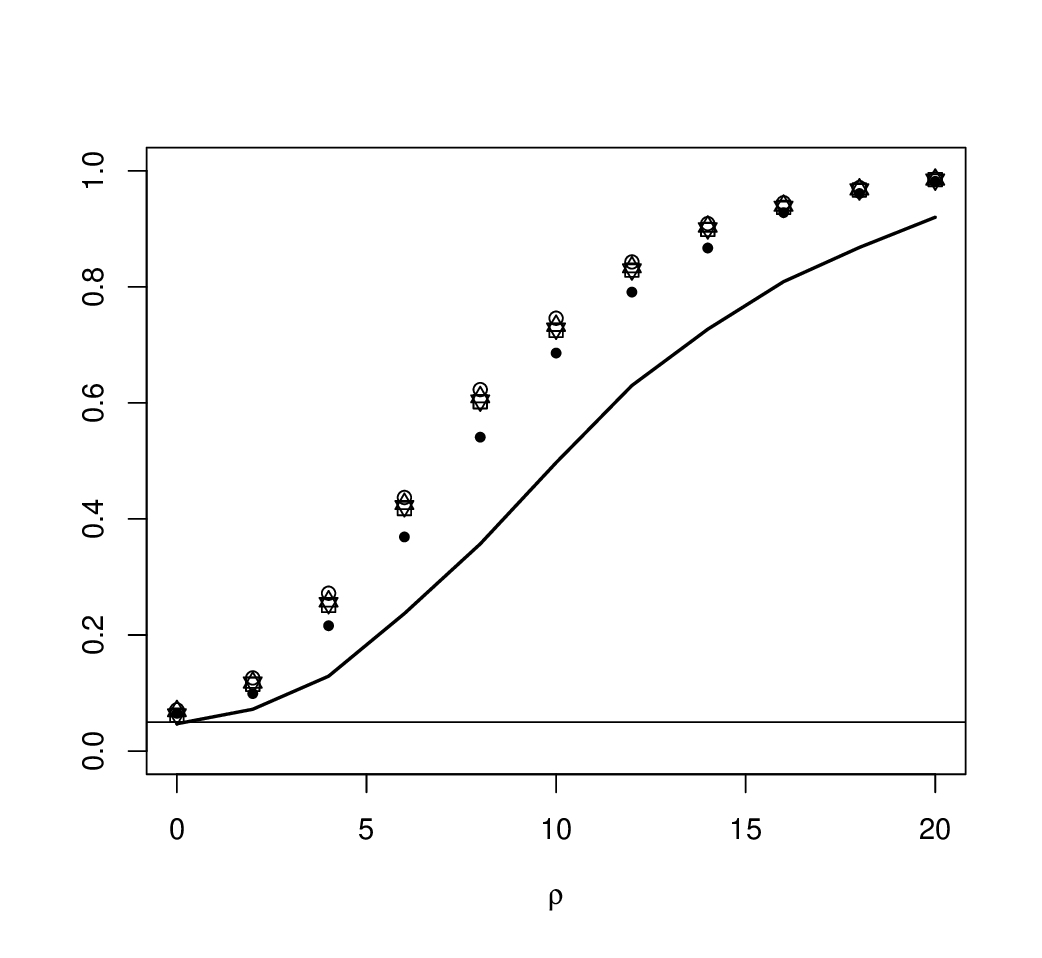} \hspace{-0.5cm}  \includegraphics[scale=0.4]{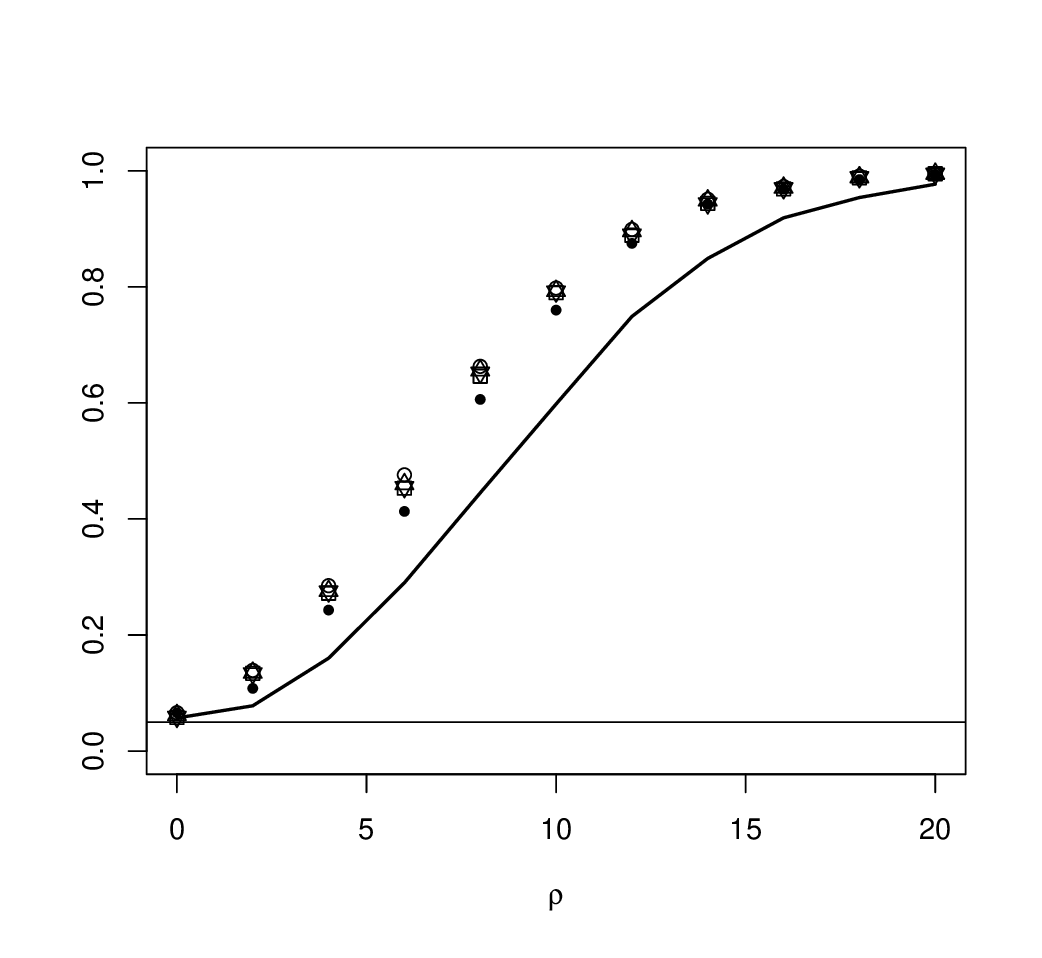}  \hspace{-0.5cm}\includegraphics[scale=0.4]{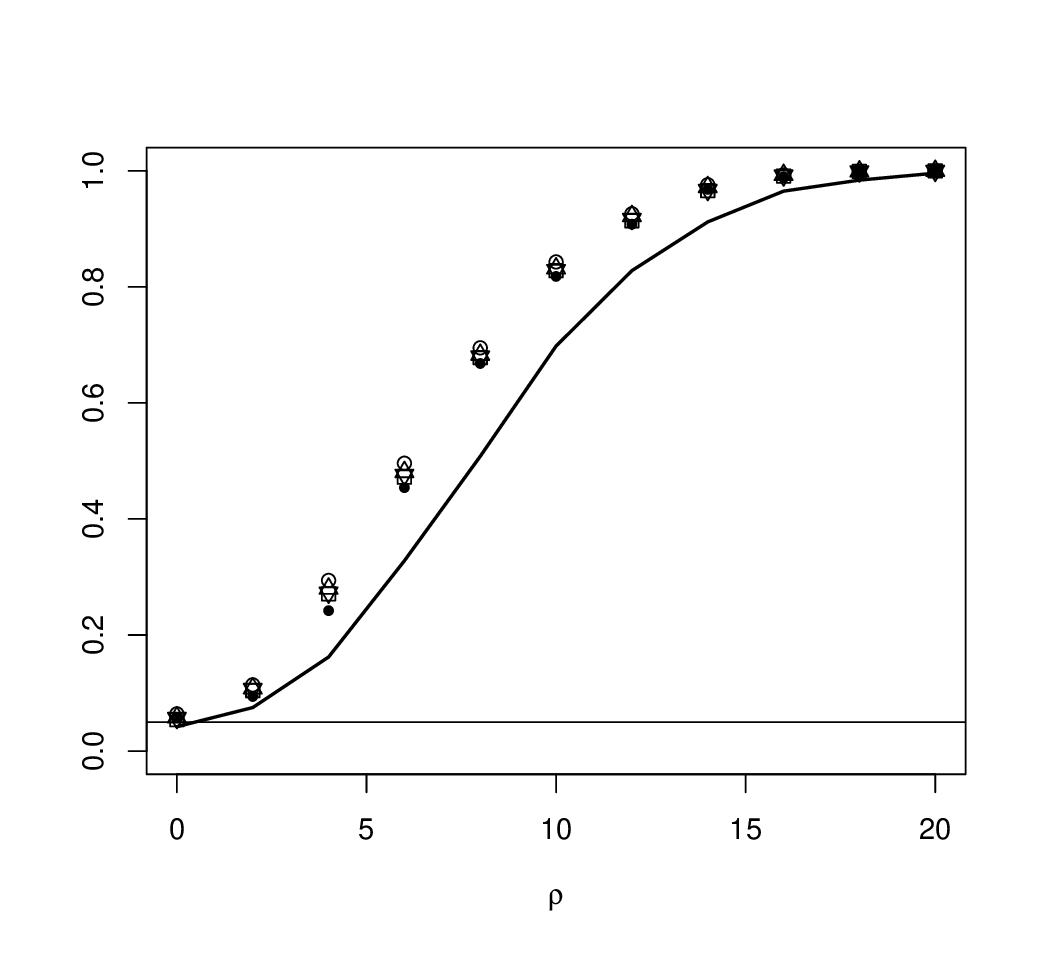} 
  \vspace{-0.5cm}
  \caption{\small \label{fig:power_k3}  Frequency of rejection when $k=3$ for the bootstrap test $\phi_{\boot,M}$, $\phi_{\boot,\gaus}$, and  $\phi_{\perm,5000}$. The solid line corresponds to $\phi_{\perm,5000}$, the filled circles to $\phi_{\boot,\gaus}$, while the circles, upper, lower triangles and the square correspond to $\phi_{\boot,M}$, with $M=3,10,20$ and $30$, respectively. }
 \end{sidewaysfigure}

 Tables \ref{tab:power_k2_perm} and \ref{tab:power_k3_perm_mod} show that the permutation test is an accurate test both for $k=2$ and $k=3$. When comparing the power of  the permutation test and the bootstrap calibration, we note that both tests lead to similar results. However, the permutation test has  a better power performance for $k=2$ when large values of $\rho$ and small sample sizes are combined. On the contrary, for $k=3$ populations a better power is attained with the bootstrap calibration.  This behaviour is clearly visualized in Figures \ref{fig:power_k2} and \ref{fig:power_k3}.

 \begin{table}[ht]
 \begin{center}
 \begin{tabular}{c|cc ||c c||c c|}
  \hline
  & \multicolumn{2}{c||}{$n_1=n_2=50$} & \multicolumn{2}{c||}{$n_1=n_2=100$} & \multicolumn{2}{c|}{$n_1=n_2=200$}\\\hline
 $\rho$  & $\phi_{\perm,1000}$ & $\phi_{\perm,5000}$ & $\phi_{\perm,1000}$ & $\phi_{\perm,5000}$ & $\phi_{\perm,1000}$ & $\phi_{\perm,5000}$\\\hline
 0 & 0.040 & 0.045 & 0.048 & 0.046 & 0.055 & 0.053 \\ 
   1 & 0.053 & 0.053 & 0.070 & 0.068 & 0.060 & 0.062 \\ 
   2 & 0.206 & 0.216 & 0.264 & 0.267 & 0.313 & 0.313 \\ 
   3 & 0.557 & 0.557 & 0.702 & 0.708 & 0.813 & 0.814 \\ 
   4 & 0.816 & 0.823 & 0.943 & 0.942 & 0.990 & 0.992 \\ 
   5 & 0.946 & 0.945 & 0.993 & 0.995 & 1.000 & 1.000 \\ 
   6 & 0.986 & 0.986 & 0.998 & 0.998 & 1.000 & 1.000 \\ 
   7 & 0.998 & 0.998 & 1.000 & 1.000 & 1.000 & 1.000 \\ 
   8 & 1.000 & 1.000 & 1.000 & 1.000 & 1.000 & 1.000 \\ 
   10 & 1.000 & 1.000 & 1.000 & 1.000 & 1.000 & 1.000 \\    \hline
 \end{tabular}
 \caption{\small \label{tab:power_k2_perm} Frequency of rejection for the permutation  test $\phi_{\perm,N_{\perm}}$   when two populations are compared and   $N_{\perm}=1000$ and $5000$ permutations are used.  The alternatives considered are $X_{1 }\sim {\cal{BW}}(0,1)$ while
 $X_{2 }\sim W_1+\delta n^{-1/4} W_2^2$, where $W_j\sim {\cal{BW}}(0,1)$ are independent of each other and $\delta_n={\rho}n^{-1/4}$ with $n= n_1+n_2$. }
 \end{center}
 \end{table}

  \begin{table}[ht]
 \centering
 \begin{tabular}{c|cc ||c c||c c|}
   \hline
   & \multicolumn{2}{c||}{$n_1=n_2=n_3=50$} & \multicolumn{2}{c||}{$n_1=n_2=n_3=100$} & \multicolumn{2}{c|}{$n_1=n_2=n_3=200$}\\\hline
  $\rho$  & $\phi_{\perm,1000}$ & $\phi_{\perm,5000}$ & $\phi_{\perm,1000}$ & $\phi_{\perm,5000}$ & $\phi_{\perm,1000}$ & $\phi_{\perm,5000}$\\ 
   \hline
  0 & 0.046 & 0.047 & 0.058 & 0.057 & 0.041 & 0.042 \\ 
   2 & 0.070 & 0.072 & 0.079 & 0.078 & 0.072 & 0.075 \\ 
   4 & 0.133 & 0.129 & 0.157 & 0.160 & 0.158 & 0.162 \\ 
   6 & 0.229 & 0.237 & 0.282 & 0.290 & 0.329 & 0.328 \\ 
   8 & 0.351 & 0.357 & 0.452 & 0.445 & 0.505 & 0.508 \\ 
   10 & 0.498 & 0.497 & 0.598 & 0.598 & 0.697 & 0.698 \\ 
   12 & 0.630 & 0.630 & 0.753 & 0.749 & 0.821 & 0.828 \\ 
   14 & 0.725 & 0.727 & 0.849 & 0.849 & 0.907 & 0.912 \\ 
   16 & 0.806 & 0.809 & 0.920 & 0.919 & 0.962 & 0.965 \\ 
   18 & 0.867 & 0.868 & 0.948 & 0.954 & 0.985 & 0.984 \\ 
   20 & 0.916 & 0.920 & 0.976 & 0.977 & 0.996 & 0.996 \\ 
   \hline
 \end{tabular}
 \caption{\small \label{tab:power_k3_perm_mod} Frequency of rejection for the permutation  test $\phi_{\perm,N_{\perm}}$  based on  $D=\|\wbGa_2 -\wbGa_1 \|_\itF^2+\|\wbGa_3 -\wbGa_1 \|_\itF^2+\|\wbGa_3 -\wbGa_2 \|_\itF^2$,   when three populations are compared,   $N_{\perm}=1000$ and $5000$ permutations are used. The observations are generated as   $X_{1}\sim {\cal{BW}}(0,1)$ and $X_{i}\sim (1+\delta_n)^{1/2}{\cal{BW}}(0,1)$, for $i=2,3$, where $\delta_n={\rho}n^{-1/2}$ with $n=\sum_{i=1}^3 n_i$. }
 \end{table}

  To help in the effective comparison of the power performance of the two tests,  we compute  the size--corrected relative exact powers $\rho_{H_1}(\phi_{\boot,M},\phi_{\perm, N_{\perm}})$ and $\rho_{H_1}(\phi_{\boot,\gaus},\phi_{\perm, N_{\perm}})$, where as mentioned above, $\phi_{\boot,M}$ stands the bootstrap calibration of $T_{k,n}$ computed using $M$ principal components, $\phi_{\boot,\gaus}$ denotes the bootstrap calibration of $T_{k,n}$ computed using the Gaussian approximation and $\phi_{\perm,N_{\perm}}$ is the permutation test computed using $N_{\perm}$ random permutations.  For two test $\phi_1$ and $\phi_2$ and an alternative $H_1$,  the size--corrected relative exact power $\rho_{H_1}(\phi_1,\phi_2)$ was defined in Morales \textsl{et al.} (2004) as
 $$
 \rho_{H_1}(\phi_1,\phi_2)=\left(\frac{D_{H_1}(\phi_1)} {D_{H_1}(\phi_2)}-1\right)\times 100\,,
 $$ 
 with $D_{H_1}(\phi)=\pi_{H_1}(\phi)-\pi_{H_0}(\phi)$, where $\pi_{H_1}(\phi)$ and $\pi_{H_0}(\phi)$ denote the power of the test $\phi$ under $H_1$ and     the null hypothesis, respectively. This measure allows to clarify the fluctuations in the powers which are more difficult to observe in Tables \ref{tab:power_k2} to \ref{tab:power_k3_perm_mod}, since large negative values of $ \rho_{H_1}(\phi_1,\phi_2)$ indicate that $\phi_2$ outperforms $\phi_1$, while large positive values show that $\phi_1$ is preferable.
 
 Tables \ref{tab:size_k2_1000} and \ref{tab:size_k2_5000} report the values of $\rho_{H_1}(\phi_{\boot, M},\phi_{\perm, N_{\perm}})$ and $\rho_{H_1}(\phi_{\boot, \gaus},\phi_{\perm, N_{\perm}})$, for two populations, when the permutation test $\phi_{\perm, N_{\perm}}$ is computed with   $N_{\perm}=1000$ and $5000$ random permutations, respectively. As expected the test obtained using the Gaussian approximation outperforms the permutation test in particular, for local alternatives close to the null hypothesis. On the other hand, the permutation test based on $1000$ permutations shows its advantage for $n_1=n_2=50$, in particular, when $\rho=1$ since the asymptotic approximation leads to some loss of power in the bootstrap test. The large negative values obtained for $N_{\perm}=1000$ are reduced when $5000$ random permutations are considered, since the empirical size is closer to the nominal one. The better performance for $\rho=1$ is also observed when $n_1=n_2=100$, while for $n_1=n_2=200$ the test defined in Section \ref{bootstrap} is much better than the permutation test. In general, for large sample sizes, the bootstrap test shows its advantage. The worst behaviour of the permutation test for large samples may be due to the fact   that the number of random permutations must be increased  with the sample size.

  \begin{sidewaystable}
 \centering
 \small
  \begin{tabular}{c|rrrr|r||rrrr|r||rrrr|r|}
   \hline
  & \multicolumn{5}{c||}{$n_1=n_2=50$} & \multicolumn{5}{c||}{$n_1=n_2=100$} & \multicolumn{5}{c|}{$n_1=n_2=200$}\\\hline
   $\rho$ &  $\phi_{\boot,3}$ &  $\phi_{\boot,10}$ &  $\phi_{\boot,20}$ &  $\phi_{\boot,30}$ & $\phi_{\boot,\gaus}$ & $\phi_{\boot,3}$ &  $\phi_{\boot,10}$ &  $\phi_{\boot,20}$ &  $\phi_{\boot,30}$ & $\phi_{\boot,\gaus}$ & $\phi_{\boot,3}$ &  $\phi_{\boot,10}$ &  $\phi_{\boot,20}$ &  $\phi_{\boot,30}$ & $\phi_{\boot,\gaus}$ \\\hline
 1 & -53.85 & -53.85 & -61.54 & -53.85 & 38.46 & -22.73 & -22.73 & -9.09 & -13.64 & -36.36 & 440.00 & 420.00 & 380.00 & 380.00 & 280.00 \\ 
   2 & 0.00 & -7.23 & -8.43 & -11.45 & 54.82 & 15.28 & 8.33 & 8.80 & 6.48 & 31.02 & 17.05 & 14.73 & 12.79 & 10.85 & 22.09 \\ 
   3 & -9.48 & -13.35 & -15.86 & -17.02 & 30.17 & -3.98 & -5.81 & -6.73 & -7.34 & 14.83 & 5.15 & 4.62 & 3.96 & 3.69 & 12.27 \\ 
   4 & -15.72 & -18.17 & -19.33 & -19.97 & 11.21 & -7.37 & -7.82 & -7.93 & -8.16 & 3.35 & 0.32 & 0.32 & 0.21 & 0.21 & 2.57 \\ 
   5 & -17.66 & -19.21 & -19.98 & -20.09 & 2.76 & -6.67 & -7.20 & -6.88 & -6.88 & 0.32 & 0.00 & 0.21 & 0.21 & 0.21 & 1.59 \\ 
   6 & -17.23 & -18.39 & -19.03 & -19.66 & 0.32 & -6.00 & -6.11 & -5.79 & -5.68 & 0.00 & 0.00 & 0.21 & 0.21 & 0.32 & 1.59 \\ 
   7 & -16.91 & -17.64 & -18.16 & -18.68 & -0.73 & -5.36 & -5.78 & -5.46 & -5.36 & -0.21 & 0.00 & 0.21 & 0.21 & 0.32 & 1.59 \\ 
   8 & -16.25 & -17.40 & -17.71 & -17.81 & -0.83 & -4.83 & -5.04 & -5.04 & -5.04 & -0.21 & 0.00 & 0.21 & 0.21 & 0.32 & 1.59 \\ 
   10 & -16.15 & -16.67 & -17.08 & -17.08 & -0.83 & -4.73 & -4.83 & -4.73 & -4.73 & -0.21 & 0.00 & 0.21 & 0.21 & 0.32 & 1.59 \\ 
    \hline
 \end{tabular}
 \caption{\label{tab:size_k2_1000} Size corrected relative exact powers, $\rho_{H_1}(\phi_{\boot, M},\phi_{\perm, N_{\perm}})$ and $\rho_{H_1}(\phi_{\boot, \gaus},\phi_{\perm, N_{\perm}})$, for the   bootstrap tests $\phi_{\boot,M}$ ($M=3,10,20$ and $30$) and $\phi_{\boot,\gaus}$ with respect to the permutation test $\phi_{\perm,N_{\perm}}$ with $N_{\perm}=1000$ random permutations, for $k=2$.} 
 
 \vskip0.2in
 \begin{tabular}{c|rrrr|r||rrrr|r||rrrr|r|}
   \hline
  & \multicolumn{5}{c||}{$n_1=n_2=50$} & \multicolumn{5}{c||}{$n_1=n_2=100$} & \multicolumn{5}{c|}{$n_1=n_2=200$}\\\hline
  $\rho$ &  $\phi_{\boot,3}$ &  $\phi_{\boot,10}$ &  $\phi_{\boot,20}$ &  $\phi_{\boot,30}$ & $\phi_{\boot,\gaus}$ & $\phi_{\boot,3}$ &  $\phi_{\boot,10}$ &  $\phi_{\boot,20}$ &  $\phi_{\boot,30}$ & $\phi_{\boot,\gaus}$ & $\phi_{\boot,3}$ &  $\phi_{\boot,10}$ &  $\phi_{\boot,20}$ &  $\phi_{\boot,30}$ & $\phi_{\boot,\gaus}$ \\  \hline
 1 & -25.00 & -25.00 & -37.50 & -25.00 & 125.00 & -22.73 & -22.73 & -9.09 & -13.64 & -36.36 & 200.00 & 188.89 & 166.67 & 166.67 & 111.11 \\ 
   2 & -2.92 & -9.94 & -11.11 & -14.04 & 50.29 & 12.67 & 5.88 & 6.33 & 4.07 & 28.05 & 16.15 & 13.85 & 11.92 & 10.00 & 21.15 \\ 
   3 & -8.59 & -12.50 & -15.04 & -16.21 & 31.45 & -5.14 & -6.95 & -7.85 & -8.46 & 13.44 & 4.73 & 4.20 & 3.55 & 3.29 & 11.83 \\ 
   4 & -15.94 & -18.38 & -19.54 & -20.18 & 10.93 & -7.48 & -7.92 & -8.04 & -8.26 & 3.24 & -0.11 & -0.11 & -0.21 & -0.21 & 2.13 \\ 
   5 & -17.11 & -18.67 & -19.44 & -19.56 & 3.44 & -7.06 & -7.59 & -7.27 & -7.27 & -0.11 & -0.21 & 0.00 & 0.00 & 0.00 & 1.37 \\ 
   6 & -16.79 & -17.96 & -18.60 & -19.23 & 0.85 & -6.20 & -6.30 & -5.99 & -5.88 & -0.21 & -0.21 & 0.00 & 0.00 & 0.11 & 1.37 \\ 
   7 & -16.47 & -17.21 & -17.73 & -18.26 & -0.21 & -5.56 & -5.97 & -5.66 & -5.56 & -0.42 & -0.21 & 0.00 & 0.00 & 0.11 & 1.37 \\ 
   8 & -15.81 & -16.96 & -17.28 & -17.38 & -0.31 & -5.03 & -5.24 & -5.24 & -5.24 & -0.42 & -0.21 & 0.00 & 0.00 & 0.11 & 1.37 \\ 
   10 & -15.71 & -16.23 & -16.65 & -16.65 & -0.31 & -4.93 & -5.03 & -4.93 & -4.93 & -0.42 & -0.21 & 0.00 & 0.00 & 0.11 & 1.37 \\ 
    \hline
 \end{tabular}
 \caption{\label{tab:size_k2_5000}Size corrected relative exact powers, $\rho_{H_1}(\phi_{\boot, M},\phi_{\perm, N_{\perm}})$ and $\rho_{H_1}(\phi_{\boot, \gaus},\phi_{\perm, N_{\perm}})$, for the   bootstrap tests $\phi_{\boot,M}$ ($M=3,10,20$ and $30$) and $\phi_{\boot,\gaus}$ with respect to the permutation test $\phi_{\perm,N_{\perm}}$ with $N_{\perm}=5000$ random permutations, for $k=2$.} 
 \end{sidewaystable}
 
 
 When considering $k=3$ populations, Tables \ref{tab:size_k3_1000_mod} and \ref{tab:size_k3_5000_mod} report  the size corrected values $\rho_{H_1}(\phi_{\boot, M},\phi_{\perm, N_{\perm}})$ and $\rho_{H_1}(\phi_{\boot, \gaus},\phi_{\perm, N_{\perm}})$   when $N_{\perm}=1000$ and $5000$, respectively. In this case, the bootstrap calibration test always outperforms the permutation test, in particular, for alternatives close to the null hypothesis. The better performance may be explained by the fact that the asymptotic behaviour of the tests and so, its bootstrap calibration, detects more easily alternatives following a proportional model than those considered in the two population case. The higher capability of $\phi_{\boot, M}$ to detect proportional local alternatives for three populations  is related to  power performance described in  Remark \ref{alt}.2. Besides, the obtained results suggest that as the number of populations increases the number of permutations needed to attain a good power performance needs also to be increased considerably, which leads to a larger computing time.
 
 
  \begin{sidewaystable}
 \centering
  \small
  \begin{tabular}{c|rrrr|r||rrrr|r||rrrr|r|}
   \hline
  & \multicolumn{5}{c||}{$n_1=n_2=n_3=50$} & \multicolumn{5}{c||}{$n_1=n_2=n_3=100$} & \multicolumn{5}{c|}{$n_1=n_2=n_3=200$}\\\hline
  $\rho$ &  $\phi_{\boot,3}$ &  $\phi_{\boot,10}$ &  $\phi_{\boot,20}$ &  $\phi_{\boot,30}$ & $\phi_{\boot,\gaus}$ & $\phi_{\boot,3}$ &  $\phi_{\boot,10}$ &  $\phi_{\boot,20}$ &  $\phi_{\boot,30}$ & $\phi_{\boot,\gaus}$ & $\phi_{\boot,3}$ &  $\phi_{\boot,10}$ &  $\phi_{\boot,20}$ &  $\phi_{\boot,30}$ & $\phi_{\boot,\gaus}$ \\\hline
 2 & 129.17 & 104.17 & 129.17 & 120.83 & 41.67 & 247.62 & 247.62 & 257.14 & 261.90 & 100.00 & 61.29 & 61.29 & 61.29 & 61.29 & 16.13 \\ 
   4 & 131.03 & 117.24 & 120.69 & 117.24 & 73.56 & 121.21 & 116.16 & 118.18 & 116.16 & 78.79 & 96.58 & 89.74 & 84.62 & 85.47 & 57.26 \\ 
   6 & 100.00 & 95.08 & 96.17 & 94.54 & 66.12 & 83.04 & 77.68 & 76.79 & 76.34 & 54.91 & 50.00 & 46.87 & 45.49 & 45.14 & 37.50 \\ 
   8 & 80.98 & 77.38 & 77.38 & 77.05 & 56.07 & 51.52 & 50.76 & 50.51 & 49.24 & 37.06 & 35.99 & 34.70 & 34.27 & 34.48 & 31.47 \\ 
   10 & 49.34 & 46.90 & 47.12 & 46.68 & 37.39 & 35.56 & 35.37 & 35.74 & 35.56 & 28.52 & 18.75 & 17.99 & 17.68 & 17.99 & 15.85 \\ 
   12 & 32.19 & 30.99 & 31.34 & 31.34 & 24.32 & 19.86 & 20.00 & 19.57 & 19.57 & 16.40 & 10.51 & 10.77 & 10.26 & 10.26 & 8.97 \\ 
   14 & 23.42 & 22.97 & 23.42 & 23.27 & 18.11 & 11.88 & 12.14 & 11.88 & 12.01 & 10.87 & 5.31 & 5.54 & 5.08 & 5.31 & 5.08 \\ 
   16 & 15.00 & 14.61 & 15.13 & 15.13 & 13.55 & 5.22 & 5.45 & 5.68 & 5.68 & 4.76 & 0.87 & 1.52 & 1.52 & 1.74 & 1.19 \\ 
   18 & 9.62 & 9.50 & 10.23 & 10.23 & 9.14 & 3.93 & 4.16 & 4.49 & 4.49 & 3.26 & -0.95 & -0.32 & -0.21 & 0.11 & -0.32 \\ 
   20 & 5.17 & 5.29 & 5.98 & 6.09 & 5.40 & 1.20 & 1.63 & 1.96 & 2.07 & 1.09 & -1.99 & -1.36 & -1.26 & -0.94 & -1.36 \\ 
    \hline
 \end{tabular}
 \caption{\small\label{tab:size_k3_1000_mod}Size corrected relative exact powers, $\rho_{H_1}(\phi_{\boot, M},\phi_{\perm, N_{\perm}})$ and $\rho_{H_1}(\phi_{\boot, \gaus},\phi_{\perm, N_{\perm}})$, for the   bootstrap tests $\phi_{\boot,M}$ ($M=3,10,20$ and $30$) and $\phi_{\boot,\gaus}$ with respect to the permutation test $\phi_{\perm,N_{\perm}}$  based on  $D=\|\wbGa_2 -\wbGa_1 \|_\itF^2+\|\wbGa_3 -\wbGa_1 \|_\itF^2+\|\wbGa_3 -\wbGa_2 \|_\itF^2$ with    $N_{\perm}=1000$ random permutations when $k=3$.} 
 
 \vskip0.2in
 
  \begin{tabular}{c|rrrr|r||rrrr|r||rrrr|r|}
   \hline
  & \multicolumn{5}{c||}{$n_1=n_2=n_3=50$} & \multicolumn{5}{c||}{$n_1=n_2=n_3=100$} & \multicolumn{5}{c|}{$n_1=n_2=n_3=200$}\\\hline
  $\rho$ &  $\phi_{\boot,3}$ &  $\phi_{\boot,10}$ &  $\phi_{\boot,20}$ &  $\phi_{\boot,30}$ & $\phi_{\boot,\gaus}$ & $\phi_{\boot,3}$ &  $\phi_{\boot,10}$ &  $\phi_{\boot,20}$ &  $\phi_{\boot,30}$ & $\phi_{\boot,\gaus}$ & $\phi_{\boot,3}$ &  $\phi_{\boot,10}$ &  $\phi_{\boot,20}$ &  $\phi_{\boot,30}$ & $\phi_{\boot,\gaus}$ \\\hline
 2 & 120.00 & 96.00 & 120.00 & 112.00 & 36.00 & 247.62 & 247.62 & 257.14 & 261.90 & 100.00 & 51.52 & 51.52 & 51.52 & 51.52 & 9.09 \\ 
   4 & 145.12 & 130.49 & 134.15 & 130.49 & 84.15 & 112.62 & 107.77 & 109.71 & 107.77 & 71.84 & 91.67 & 85.00 & 80.00 & 80.83 & 53.33 \\ 
   6 & 92.63 & 87.89 & 88.95 & 87.37 & 60.00 & 75.97 & 70.82 & 69.96 & 69.53 & 48.93 & 51.05 & 47.90 & 46.50 & 46.15 & 38.46 \\ 
   8 & 78.06 & 74.52 & 74.52 & 74.19 & 53.55 & 53.87 & 53.09 & 52.84 & 51.55 & 39.18 & 35.41 & 34.12 & 33.69 & 33.91 & 30.90 \\ 
   10 & 50.00 & 47.56 & 47.78 & 47.33 & 38.00 & 35.30 & 35.12 & 35.49 & 35.30 & 28.28 & 18.75 & 17.99 & 17.68 & 17.99 & 15.85 \\ 
   12 & 32.42 & 31.22 & 31.56 & 31.56 & 24.53 & 20.38 & 20.52 & 20.09 & 20.09 & 16.91 & 9.67 & 9.92 & 9.41 & 9.41 & 8.14 \\ 
   14 & 23.24 & 22.79 & 23.24 & 23.09 & 17.94 & 11.74 & 11.99 & 11.74 & 11.87 & 10.73 & 4.83 & 5.06 & 4.60 & 4.83 & 4.60 \\ 
   16 & 14.70 & 14.30 & 14.83 & 14.83 & 13.25 & 5.22 & 5.45 & 5.68 & 5.68 & 4.76 & 0.65 & 1.30 & 1.30 & 1.52 & 0.98 \\ 
   18 & 9.62 & 9.50 & 10.23 & 10.23 & 9.14 & 3.12 & 3.34 & 3.68 & 3.68 & 2.45 & -0.74 & -0.11 & 0.00 & 0.32 & -0.11 \\ 
   20 & 4.81 & 4.93 & 5.61 & 5.73 & 5.04 & 0.98 & 1.41 & 1.74 & 1.85 & 0.87 & -1.89 & -1.26 & -1.15 & -0.84 & -1.26 \\ 
    \hline
 \end{tabular}
 \caption{\small\label{tab:size_k3_5000_mod}Size corrected relative exact powers, $\rho_{H_1}(\phi_{\boot, M},\phi_{\perm, N_{\perm}})$ and $\rho_{H_1}(\phi_{\boot, \gaus},\phi_{\perm, N_{\perm}})$, for the   bootstrap tests $\phi_{\boot,M}$ ($M=3,10,20$ and $30$) and $\phi_{\boot,\gaus}$ with respect to the permutation test $\phi_{\perm,N_{\perm}}$  based on  $D=\|\wbGa_2 -\wbGa_1 \|_\itF^2+\|\wbGa_3 -\wbGa_1 \|_\itF^2+\|\wbGa_3 -\wbGa_2 \|_\itF^2$ with    $N_{\perm}=5000$ random permutations when $k=3$.} 
 \end{sidewaystable}
 
 Although a formal computational  analysis of the different test statistics is beyond the scope of this paper, we tested the speed
 of our R codes using   an Intel i7-2600K CPU (3.4GHz) machine running Windows 7.  Table \ref{tab:tiempos} report the average time in CPU seconds of the different test procedures computed over 10 random samples generated   as in the simulation settings under $H_0$ and for the sample sizes $n_i$ considered above.  The obtained results show that the computing  time increases linearly as  the number of permutations increase and in all situations $\phi_{\perm, N_{\perm}}$ is much more time expensive than $\phi_{\boot,M}$. On the other hand, as expected, the number $M$ of principal components used increases considerably the computation time. However,   the computing time of  $\phi_{\boot,M}$ is quite stable along sample sizes, for a fixed number of populations and a fixed $M$. The Gaussian approximation takes almost the same computing time in all the considered situations and shows a larger average time than $\phi_{\boot,M}$, except when $M=30$ and $k=3$, in which they both give similar average timings.

  \begin{table}[ht]
 \centering
 \begin{tabular}{c|cc c||c c  c|}
   \hline
  & \multicolumn{3}{c|}{$k=2$} & \multicolumn{3}{c|}{$k=3$}\\ \hline
  &   {$n_i=50$} &  {$n_i=100$} &  {$n_i=200$} & {$n_i=50$} &  {$n_i=100$} &  {$n_i=200$}\\\hline
  $\phi_{\boot,3}$ & 0.053 & 0.059 & 0.090 & 0.055 & 0.072 & 0.114 \\ 
   $\phi_{\boot,10}$ & 0.125 & 0.120 & 0.151 & 0.164 & 0.173 & 0.215 \\ 
   $\phi_{\boot,20}$& 0.334 & 0.309 & 0.367 & 0.693 & 0.693 & 0.828 \\ 
  $\phi_{\boot,30}$ & 0.867 & 0.906 & 1.069 & 3.510 & 3.580 & 3.822 \\ \hline
  $\phi_{\boot,\gaus}$ & 3.424 & 3.363 & 3.379 & 3.317 &    3.315 & 4.413 \\ \hline
 $\phi_{\perm,1000} $ & 1.176 & 1.930 & 3.457 & 3.264 & 5.045 & 9.276 \\ 
 $\phi_{\perm,5000} $& 5.831 & 9.493 & 17.825 & 15.544 & 25.957 & 47.575 \\
    \hline
 \end{tabular}
 \caption{\small \label{tab:tiempos} Average timing (in seconds) of the test procedures.}
 \end{table}

 From the obtained results, we see that our procedure is, in terms of level and power behaviour,  a good competitor for the permutation test introduced for two populations in Pigoli \textsl{et al.} (2014). On the other hand, when $k=3$ it has a better detection capability with a much lower computing time.  Besides, our method has the advantage of allowing to develop a theory regarding its asymptotic behaviour as described  in Sections \ref{alt} and \ref{bootstrap}.

 \section{Conclusions}
 
 In this paper, we have studied a procedure to test equality among several populations covariance operators. The test statistic is based on the Hilbert--Schmidt distance between consistent estimators of  $\bGa_i$ and $\bGa_1$, for $2\le i\le k$.  The analysis of the asymptotic distribution of the test statistic reveals that the testing procedure is consistent against local
 alternatives converging to the null hypothesis at rate $n^{-1/2}$, when the sample covariance operators are used. These results also hold for the smoothed covariance operators defined in Boente and Fraiman (2000), under mild conditions. The asymptotic null behaviour obtained motivate the use of   bootstrap methods, since it depends on the eigenvalues of an unknown operator. For that reason,
 we also provide a general bootstrap calibration method whose validity is derived. Our numerical studies  have shown that the bootstrap calibration has  a good practical behaviour and is a good competitor for the permutation test defined in Pigoli et al. (2014) for two populations and the considered alternatives.  On the other hand, when $k=3$ and for proportional alternatives, it has shown a better detection capability. Another advantages of the bootstrap test over the permutation test is its lower computing time, for the sample sizes considered.

 \bigskip
 
 \noi\textbf{\small Acknowledgements.} {\small This research was partially supported by Grants  \textsc{pip}
 112-201101-00339 and 112-201101-00742 from \textsc{conicet}, \textsc{pict} 2014-0351 and 2012-1641 from
 \textsc{anpcyt} and  20020130100279\textsc{ba} and 20020120200244\textsc{ba}  from the Universidad de
 Buenos Aires  at Buenos Aires, Argentina.}
 
  \setcounter{section}{0}
 \renewcommand{\thesection}{A}

 \section{Appendix}{\label{ape}}
 \setcounter{equation}{0}
 \renewcommand{\theequation}{A.\arabic{equation}}
 
 \noi \textbf{Proof of Theorem \ref{stest}.1.} 
 Denote as $\itF^{k}=\itF\times \dots \times \itF$ the $k-$th dimensional product space of identical copies of $\itF$ and  consider the process $\bV_{k,n}=\left(\sqrt n (\wtbGa_1-\bGa_1), \dots, \sqrt n (\wtbGa_k-\bGa_k)\right)\trasp$. Using that $\sqrt n_i\left({\wtbGa_i}-\bGa_i\right)\convdist \bU_i$, the independence of the estimated operators  and the fact that  $n_i/n\to \tau_i\in (0,1)$, we get that
 $\bV_{k,n}\convdist \bV=(\bV_1,\cdots,\bV_k)\trasp$,
 where $\bV_i=\tau_i^{-1/2}\bU_i$ are independent random processes of $\itF$  with covariance operators ${\tau_i}^{-1}\bUpsi_i$.  Hence, $\bV_{k,n}$ converges in distribution to a zero mean Gaussian random element $\bV=(\bV_1,\cdots,\bV_k)\trasp\in \itF^{k}$   with covariance operator $\wtbUps=\diag\left( {\tau_1}^{-1}\bUpsi_1,\dots, {\tau_k}^{-1}\bUpsi_k\right)$.
 
 Let $A:\itF^{k}\to \itF^{k-1}$ be the linear operator defined as $A(V_1,\cdots,V_k)=(V_2-V_1,\cdots,V_k-V_1) $  and denote as $A^*:\itF^{k-1}\to \itF^{k}$ its adjoint operator.
 The continuous map theorem guarantees that $ A\bV_{k,n}\convdist \bW$, where $\bW=(W_1,\dots, W_{k-1})\trasp= A \bV$ is a zero mean Gaussian random element of $\itF^{k-1}$ with covariance operator $\bUpsi_{\weigh}=A\wtbUps\, A^*$. Moreover, we also obtain that $n \sum_{j=2}^k\|(\wtbGa_j-\bGa_j)-(\wtbGa_1-\bGa_1)\|_\itF^2\convdist \sum_{j=1}^{k-1} \|W_j\|_\itF^2=\|\bW\|_{\itF^{k-1}}^2$. Let $\upsilon_\ell\in \itF^{k-1}$ be the orthonormal eigenfunctions of $\bUpsi_{\weigh}$ related to the eigenvalues $\theta_\ell$ ordered in decreasing order. Since $\bW$ is a zero mean Gaussian random element  of $\itF^{k-1}$  with covariance operator $\bUpsi_{\weigh}$,   $\bW$ can be
 written as $\sum_{\ell\,\ge 1} \theta_\ell^{1/2} Z_\ell\,\upsilon_\ell$ where $Z_\ell$
 are i.i.d. random variables such that $Z_\ell\sim N(0,1)$. Hence, $\|\bW\|_{\itF^{k-1}}^2=\sum_{\ell\,\ge 1} \theta_\ell Z_\ell^2$, which leads to the desired result.
 
  It only remains to show (\ref{upsiweight}).
 Straightforward  calculations allow to show that the adjoint operator   $A^*:\itF^{k-1}\to \itF^k$ is
 given by
  $
  A^*(w_1,\dots,w_{k-1})\;=\;(-\sum_{i=1}^{k-1}w_i,w_1,\dots,w_{k-1})$. Hence,
 as $\bU_1,\cdots,\bU_k$ are independent, we obtain that
 \begin{eqnarray*}
  \bUpsi_{\weigh}(w_1,\dots,w_{k-1})&=&(A\wtbUps\, A^*)(w_1,\dots,w_{k-1})\\
  &=& \left(\frac{1}{\tau_2}\bUpsi_2(w_1)+\frac{1}{\tau_1}\bUpsi_1 \left(\sum_{i=1}^{k-1}w_i\right),\dots,\frac{1}{\tau_{k}}\bUpsi_{k} (w_{k-1})+\frac{1}{\tau_1}\bUpsi_1 \left(\sum_{i=1}^{k-1}w_i\right)\right) \;,
  \end{eqnarray*}
 concluding the proof. \square
 
 \vskip0.1in

 \noi \textbf{Proof of Corollary \ref{stest}.1.}  Consider the process $ \bU_{i,n_i}=\sqrt n_i (\wbGa_i-\bGa_i) $. The independence of the samples and among populations together with  the results stated in Dauxois \textsl{et al.} (1982), allow to show that $ \bU_{i,n_i}$ are independent and converge in
 distribution to independent  zero mean Gaussian random elements $\bU_i$ of $\itF$ with covariance operator $\bUpsi_i$ defined in (\ref{varasintgamai}). The result follows now from Theorem \ref{stest}.1. \square
 
 \vskip0.1in
 \noi \textbf{Proof of Theorem \ref{alt}.1.}
 Using that $n_i/n\to \tau_i$, we get immediately  that $\sqrt {n} \left(\wtbGa_i-\bGa_1\right)\convdist \bDelta_i+ (1/\sqrt {\tau_i})\bU_i$ where $\bU_i$ is a zero mean Gaussian random element with covariance operator $\bUpsi_i$ and  for $i=1$,  $\bDelta_1=\bO$ is the null operator.   The fact that the estimators are independent implies that $\bU_i$   can be chosen to be independent so, as in the proof of Theorem \ref{stest}.1, we have that $\bV_{k,n}=\left(\sqrt n (\wtbGa_1-\bGa_1), \dots, \sqrt n (\wtbGa_k-\bGa_1)\right)\trasp\convdist \bV=(\bV_1,\cdots,\bV_k)\trasp$,
 where $\bV_i= \bDelta_i+ (1/\sqrt {\tau_i})\bU_i$ are independent random processes of $\itF$  with mean $ \bDelta_i$ and covariance operators ${\tau_i}^{-1}\bUpsi_i$.  Hence, $\bV_{k,n}$ converges in distribution to a   Gaussian random element $\bV=(\bV_1,\cdots,\bV_k)\trasp\in \itF^{k}$   with mean $\bDelta=(\bDelta_1,\dots, \bDelta_k)\trasp$ and covariance operator $\wtbUps=\diag\left( {\tau_1}^{-1}\bUpsi_1,\dots, {\tau_k}^{-1}\bUpsi_k\right)$.
 
 As in the proof of Theorem \ref{stest}.1,   define $A:\itF^{k}\to \itF^{k-1}$ as the linear operator  $A(V_1,\cdots,V_k)=(V_2-V_1,\cdots,V_k-V_1) $. Then,  $ A\bV_{k,n}\convdist \bW$, where $\bW=(W_1,\dots, W_{k-1})\trasp= A \bV$ is a   Gaussian random element of $\itF^{k-1}$ with mean $A\bDelta$ and covariance operator $A\wtbUps\, A^*$. Note that $A\bDelta=(\bDelta_2,\dots, \bDelta_k)=\bDelta^{(k-1)}$, since $\bDelta_1$ is the null operator. Moreover, from the proof of Theorem \ref{stest} we get that $A\wtbUps\, A^*=\bUpsi_{\weigh}$.  Let $\upsilon_\ell\in \itF^{k-1}$ be the orthonormal eigenfunctions of $\bUpsi_{\weigh}$ related to the eigenvalues $\theta_\ell$ ordered in decreasing order. Since $\bW- \bDelta^{(k-1)}$ is a zero mean Gaussian random element  of $\itF^{k-1}$  with covariance operator $\bUpsi_{\weigh}$,   $\bW- \bDelta^{(k-1)}$ can be
 written as $\sum_{\ell\,\ge 1} \theta_\ell^{1/2} Z_\ell\,\upsilon_\ell$ where $Z_\ell$
 are i.i.d. random variables such that $Z_\ell\sim N(0,1)$. On the other hand, we have the expansion $\bDelta^{(k-1)}=\sum_{\ell\,\ge 1} \eta_\ell\,\upsilon_\ell$, so that $\bW= \sum_{\ell\,\ge 1} \left(\eta_\ell+\theta_\ell^{1/2} Z_\ell\right)\upsilon_\ell$ and  $\|\bW\|_{\itF^{k-1}}^2=\sum_{\ell\,\ge 1} \left(\eta_\ell+\theta_\ell^{1/2} Z_\ell\right)^2=\sum_{\ell\ge 1}
 \theta_\ell\left(\eta_\ell\,\theta_\ell^{-1/2}+ Z_\ell\right)^2$, which concludes the proof since  $T_{k,n}=n \sum_{j=2}^k\|(\wbGa_j-\bGa_1)-(\wbGa_1-\bGa_1)\|_\itF^2\convdist \sum_{j=1}^{k-1} \|W_j\|_\itF^2=\|\bW\|_{\itF^{k-1}}^2$. \square
 
 \vskip0.1in
 \noi \textbf{Proof of Proposition \ref{alt}.1.}  
 The results in  Dauxois \textsl{et al.} (1982) entail that $\sqrt{n_1} \left(\wbGa_{1}-\bGa_1\right)\convdist \bU_1$, where $\bU_1$ a zero mean Gaussian random element with covariance operator $\bUpsi_1$ so, we only have to prove the result for $i\ge 2$.
 Note that 
 \begin{eqnarray*}
 \sqrt{n_i}(\wbGa_i-\bGa_1)\!\!&=&\sqrt{n_i}\left(\frac1{n_i}\sum_{j=1}^{n_i}(X_{i,j}-
 \overline{X}_{i})\otimes(X_{i,j}-\overline{X}_{i})-\bGa_1\right)\\
 &=&\sqrt{n_i}(\wtbGa_i-\bGa_1)+{n^{-1/4}}{\sqrt{n_i}}\;\wbGa_{i,WR}+{n^{-1/4}}{\sqrt{n_i}}\;\wbGa_{i,RW}+{n^{-1/2}}\sqrt{n_i}\;\wbDelta_i
 \end{eqnarray*}
 where 
 \begin{eqnarray*}
 \wtbGa_i = \frac1{n_i}\sum_{j=1}^{n_i}(W_{i,j}-\overline{W}_{i})\otimes(W_{i,j}-\overline{W}_{i}),&\qquad & \wbDelta_i = \frac1{n_i}\sum_{j=1}^{n_i}(R_{i,j}-\overline{R}_{i})\otimes(R_{i,j}-\overline{R}_{i})\,,\\
 \wbGa_{i,WR} = \frac1{n_i}\sum_{j=1}^{n_i}(W_{i,j}-\overline{W}_{i})\otimes(R_{i,j}-\overline{R}_{i})  &\mbox{ and } &  
 \wbGa_{i,RW} = \frac1{n_i}\sum_{j=1}^{n_i}(R_{i,j}-\overline{R}_{i})\otimes(W_{i,j}-\overline{W}_{i})\,.
 \end{eqnarray*}
 Using that $W_{i,j}\sim W_i$ and that the covariance operator of $W_i$ is $\bGa_1$, from the results in  Dauxois \textsl{et al.} (1982) we get that 
 $\sqrt{n_i} \left(\wtbGa_{i}-\bGa_1\right)\convdist \bU_i$, where $\bU_i$ a zero mean Gaussian random element with covariance operator $\bUpsi_i$ given in (\ref{UPSIi1}).

  Note that $\wbGa_{i,WR}$ and $\wbGa_{i,RW}$ are estimators  of the cross covariance operators
  $\bGa_{i,WR}=\esp\left\{(W_{i}-\esp {W}_{i})\otimes\right.$ $\left.(R_{i}- \esp {R}_{i})\right\}$ and 
  $\bGa_{i,RW}=\esp\left\{(R_{i}-\esp{R}_{i})\otimes(W_{i}-\esp{W}_{i})\right\}$, respectively. The independence between  $W_{i}$ and  $R_{i}$  entails that  $\bGa_{i,WR}$ is the null operator, which implies that  $\sqrt{n_i}\;\wbGa_{i,WR}$ is bounded in probability, so that $n^{-1/4}\sqrt{n_i}\;\wbGa_{i,WR}\convprob 0$. Similarly, we obtain that $n^{-1/4}\sqrt{n_i}\;\wbGa_{i,RW}\convprob 0$.
 
 Finally, using the law of large numbers we have that  $\wbDelta_i$, the empirical covariance operator of $R_{i}$,
 converges in probability to $\bDelta_i$, so ${n^{-1/2}}\sqrt{n_i}\;\wbDelta_i\convprob \tau_i^{1/2}\bDelta_i$, concluding the proof of a). \square
 
 \vskip0.1in
   
 \noi \textbf{Proof of Proposition \ref{alt}.2.}  
 As in the proof of Proposition \ref{alt}.1, we only have to prove the result for $i\ge 2$.
 Using the Karhunen--Lo\'eve representation, we can write
 \begin{eqnarray*}
  X_{1,j}&=&\mu_1+\;\sum_{\ell=1}^\infty \lambda_\ell^{\frac
  12}\,\,f_{1\ell j}\,\,\phi_\ell\;,\;1\le j\le n_1\\
   X_{i,j}&=&\mu_i+\;\sum_{\ell=1}^\infty \lambda_\ell^{\frac 12}\,\left(1+\frac{\Delta_{i,\ell}}
   {\sqrt{n}}\right)^{\frac 12}\,f_{i\ell j}\,\phi_\ell\,, \quad 1\le j\le n_i\,,\quad 2\le i\le k\;.
  \end{eqnarray*}
  where $f_{i\ell j}\sim f_{i\ell}$ in (\ref{fcpc}).
 For $1\le j\le n_i$, let $Z_{i,j}=\mu_i+\;\sum_{\ell=1}^\infty \lambda_\ell^{\frac 12}\,f_{i\ell j}\,\phi_\ell=\mu_i+Z_{0,i,j}$. Denote as
 $$V_{i,j}= X_{i,j}-Z_{i,j}=\sum_{\ell=1}^\infty \lambda_\ell^{\frac 12}\left[\left(1+\frac{\Delta_{i,\ell}}{\sqrt{n}}\right)^{\frac 12}-1\right]\,f_{i\ell j}\,\phi_\ell\;.$$
 Define the following operators that will be used in the sequel $\wtbGa_i= (1/{n_i}) \sum_{j=1}^{n_i} \left(X_{i,j}-\mu_i\right)\otimes \left(X_{i,j}-\mu_i\right)$ , $\wbGa_{Z_0}= (1/{n_i}) \sum_{j=1}^{n_i} Z_{0,i,j}\otimes Z_{0,i,j}$,  $\wbGa_V = (1/{n_i}) \sum_{j=1}^{n_i} V_{i,j}\otimes V_{i,j}$ and finally, $ \wtbA = (1/{n_i}) \sum_{j=1}^{n_i} (Z_{0,i,j}\otimes V_{i,j}+V_{i,j}\otimes Z_{0,i,j})$, where we avoid the index $i$ for the sake of simplicity. Using that $X_{i,j}-\mu_i=Z_{0, i, j} +V_{i, j}$, we obtain the following expansion $\wtbGa_i=\wbGa_{Z_0}+ \wbGa_V+\wtbA$.
 
 The proof will be carried out in several  steps, by showing that
 \begin{eqnarray}
 \sqrt{n_i}(\wbGa_i-
 \wtbGa_i)&=&o_{\prob}(1)\label{reemplazo}\\
 \sqrt{n_i} \,\wbGa_V &=&o_{\prob}(1)\label{convgav}\\
 \sqrt{n_i}\, \wtbA&\convprob & \tau_i^{\frac 12} \bDelta_i \label{convergeA}\\
 \sqrt{n_i}\,(\wbGa_{Z_0}-\bGa_1)& \convdist& \bU_i\;,
 \label{convgaz0}
 \end{eqnarray}
 where $\bU_i$ is a zero mean Gaussian random element with covariance
 operator $\bUpsi_i$. Using that, for all $2\le i \le k$,  the covariance operator of $Z_{0,i,j}$
 is $\bGa_1$,  (\ref{convgaz0}) follows from Dauxois \textsl{et al.} (1982).

 We will derive (\ref{reemplazo}). Note that  $\wbGa_i-
 \wtbGa_i=-\left(\overline{X}_{i}-\mu_i\right)\otimes
 \left(\overline{X}_{i}-\mu_i\right)$. Then, it is enough  to prove
 that $\sqrt{n_i} \left(\overline{X}_{i}-\mu_i\right)=\sqrt{n_i}
 \left(\overline{Z}_{0,i} +\overline{V}_i \right)=O_\prob(1)$, with $\overline{Z}_{0,i}=(1/{n_i}) \sum_{j=1}^{n_i} Z_{0,i,j}$ and $\overline{V}_i=(1/{n_i}) \sum_{j=1}^{n_i} V_{i,j}$.
 
 By the
 central limit theorem in Hilbert spaces, we get that
 $\sqrt{n_i}\,\overline{Z}_{0,i}$ converges in distribution, which entails that the process is
 tight, i.e., $\sqrt{n_i}\,\overline{Z}_{0,i}=O_{\prob}(1)$.

 Note that
 \begin{equation}
 \left(1+\frac{\Delta_{i,\ell}}{\sqrt{n}}\right)^{\frac 12}-1= \frac 1{\sqrt{n}}\, \frac{\Delta_{i,\ell}}{\left(1+\frac{\Delta_{i,\ell}}{\sqrt{n}}\right)^{\frac 12}+1}= a_{i,\ell,n}\frac {\Delta_{i,\ell}}{\sqrt{n}}
 \label{otra}
 \end{equation}
 where $0\le a_{i,\ell,n}\le 1$.
 
 To derive that $\sqrt{n_i}\,\overline{V}_i=O_\prob(1)$, we will further show that
 $\sqrt{n_i}\,\overline{V}_i=o_\prob(1)$. To do so, note that $
 \esp \,\| \overline{V}_i\|^2 =({1}/{n_i})\sum_{\ell=1}^\infty \lambda_\ell\left[\left(1+({\Delta_{i,\ell}}/{\sqrt{n}})\right)^{\frac 12}-1\right]^2$.
 Using (\ref{otra}), we get that $\esp(\|\sqrt{n_i} \,\overline{V}_i\|^2)\leq ({1}/{n}) \sum_{\ell=1}^\infty \lambda_\ell  \Delta_{i,\ell}^2 $, concluding the proof of (\ref{reemplazo}).

 To obtain (\ref{convgav}), note that  (\ref{otra}) entails that $V_{i,j}\otimes V_{i,j}
 =(1/{n})\sum_{\ell, s} \lambda_\ell^{\frac 12}\lambda_s^{\frac 12} a_{i,\ell,n}     a_{i,s,n}   \Delta_{i,s} \Delta_{i,\ell}   f_{i \ell j} f_{i s j}\phi_\ell \otimes\phi_s$,
 so if we denote  as $U_{\ell s}=({1}/{n_i})\sum_{j=1}^{n_i}f_{i \ell j}\, \,f_{i
 sj}$, we get that
 $$
 \wbGa_V= \frac 1{n} \sum_{\ell, s} \lambda_\ell^{\frac 12}\lambda_s^{\frac 12}a_{i,\ell,n}   \, a_{i,s,n}   \Delta_{i,s}   \Delta_{i,\ell} \, U_{\ell s}\,\;\phi_\ell \otimes\phi_s\,.$$
  Note that $f_{i \ell j}\sim f_{i \ell }$ and recall that  $\esp(f_{i \ell }f_{i s})=\delta_{\ell s}$ where $\delta_{\ell s}=1$ if $\ell =s$ and $0$ otherwise. Hence, we have that $\esp (U_{\ell s})= \delta_{\ell s}$ which implies that
 \begin{eqnarray}
 \esp(U^2_{\ell s})&=&\var (U_{\ell s})+ \esp^2 (U_{\ell s})=\frac{1}{n_i}\var(f_{i\ell }f_{is}) + \delta_{\ell s}\nonumber \\
 & \leq & \frac{1}{n_i}\esp(f^2_{i\ell }f^2_{is}) + \delta_{\ell s} \leq \frac{1}{n_i}\sigma_{4,i,\ell} \sigma_{4,i,s} + \delta_{\ell s}\;,
 \label{espUls}
 \end{eqnarray}
 where the last bound follows from the Cauchy--Schwartz inequality and the fact that  $\sigma^2_{4,i,s}=\esp(f^4_{is})$. Hence, using (\ref{espUls}) and the fact that $0\le a_{i,\ell,n}\le 1$, we obtain the bound
 \begin{eqnarray*}
 \esp(n_i\|\wbGa_V\|_\itF^2)&\le &\frac{n_i}{n^2}\sum_{\ell, s} \lambda_\ell\lambda_s     \Delta_{i,s}^2  \Delta_{i,\ell}^2\, \esp(U^2_{\ell s}) \le  \frac{n_i}{n^2}\sum_{\ell, s} \lambda_\ell\lambda_s \,\Delta_{i,\ell}^2 \Delta_{i,s}^2\left(\frac{1}{n_i}\sigma_{4,i,\ell}\sigma_{4,i,s}+ \delta_{\ell s}\right)\\
 & & =\frac{1}{n^2}\left(\sum_{\ell} \lambda_\ell \,{\Delta_{i,\ell}^2} {\sigma_{4,i,\ell}}\right)^2+\frac{1}{n} \sum_{\ell} \lambda_\ell^2\,\Delta_{i,\ell}^4\,.
 \end{eqnarray*}
 Therefore, from the fact that $ \sum_{\ell} \lambda_\ell^2\,\Delta_{i,\ell}^4\le  \left(\sum_{\ell} \lambda_\ell \,\Delta_{i,\ell}^2\right)^2<\infty$ we get that $\esp({n_i}\| \wbGa_V\|_\itF^2)\to 0$, concluding the proof of (\ref{convgav}).
 
 Finally, to derive (\ref{convergeA}) note that
 $$\wtbA=\frac{1}{\sqrt{n}} \sum_{\ell, s} \lambda_\ell^{\frac 12}\lambda_s^{\frac 12}a_{i,s,n}       \Delta_{i,s}  \, U_{\ell s}\,\left(\phi_\ell\;\otimes\phi_s+\phi_s\;\otimes\phi_\ell\right)=\wtbA_1+\wtbA_2$$
 where $U_{\ell s}=({1}/{n_i})\sum_{j=1}^{n_i}f_{i \ell j}\, \,f_{isj}$, as above. We will show that $\sqrt{n_i}\left(\wtbA_j-\esp(\wtbA_j)\right)\convprob 0$, for $j=1,2$ which entails that $\sqrt{n_i}\left(\wtbA-\esp(\wtbA)\right)\convprob 0$. We will only prove that $\sqrt{n_i}\left(\wtbA_1-\esp(\wtbA_1)\right)\convprob 0$, since the other one follows similarly.  Note that from (\ref{espUls}), we get that $\var(U_{\ell s})\le ({1}/{n_i})\sigma_{4,i,\ell} \sigma_{4,i,s} $ which together with the fact that $0\le a_{i,\ell,n}\le 1$ leads to
 \begin{eqnarray*}
 \esp(n_i\|\wtbA_1-\esp(\wtbA_1)\|_\itF^2)&=&  \frac{n_i}{n} \sum_{\ell, s} \lambda_\ell \lambda_s a_{i,s,n}^2       \Delta_{i,s}^2 \var(U_{\ell s})\\
 &\le &\frac{1}{n} \sum_{\ell, s} \lambda_\ell \lambda_s        \Delta_{i,s}^2   \sigma_{4,i,\ell} \sigma_{4,i,s}  =
   \frac{1}{\sqrt n}\left(\sum_{\ell}\lambda_\ell\sigma_{4,i,\ell}\right)\left(\sum_{\ell}\lambda_\ell\sigma_{4,i,\ell}\Delta_{i,\ell}^2\right)
 \end{eqnarray*}
  so that  $\sqrt{n_i}\left(\wtbA_1-\esp(\wtbA_1)\right)\convprob 0$, as desired. Besides, using that  $\esp(U_{\ell s})=\delta_{\ell s}$, we get that
 \begin{eqnarray*}
 \esp(\sqrt{n_i}\, \wtbA) &=& \frac{2\, \sqrt{n_i}}{\sqrt{n}} \sum_{\ell} \lambda_\ell\, a_{i,\ell,n}          \Delta_{i,\ell} \; \phi_\ell\otimes \phi_\ell \to \tau_i^{\frac 12} \sum_{\ell=1}^\infty \lambda_\ell {\Delta_{i,\ell}}\phi_\ell\otimes \phi_\ell= \tau_i^{\frac 12} \bDelta_i
 \end{eqnarray*}
 where  we have used that $a_{i,\ell,n} \to 1/2$, as $n\to \infty$ and  $ \sum_{\ell=1}^\infty \lambda_\ell |{\Delta_{i,\ell}}|<\infty$. This concludes   the proof of (\ref{convergeA}). The proof of Proposition \ref{alt}.2a) follows now combining (\ref{reemplazo}) to (\ref{convgaz0}). \square

 \vskip0.1in
 
 \noi \textbf{Proof of Theorem \ref{bootstrap}.1.}
 Recall that $\wtX_n=(X_{1,1},\cdots,X_{1,n_1},\dots, X_{k,1},\cdots,X_{k,n_k})$. Let
  $\wtZ_n=(Z_1,\cdots,Z_{q_n})$ and $\wtZ =\{Z_\ell\}_{\ell\,\ge 1}$ with $Z_i\sim N(0,1)$ independent. Define
 $\widehat{\itU}_n(\wtX_n,\wtZ_n)=\sum_{\ell=1}^{q_n}\wtheta_\ell Z_\ell^2$, ${\itU}_n(\wtZ_n)=\sum_{\ell=1}^{q_n} \theta_\ell  Z_\ell^2$ and ${\itU}(\wtZ)=\sum_{\ell=1}^{\infty} \theta_\ell  Z_\ell^2$.
 
 First note that, for any $\ell$,   $|\wtheta_\ell -\theta_\ell | \leq \|\wbUps_{\weigh}-\bUpsi_{\weigh}\|_{\itG^{k-1}}$
 (see, for instance, Kato, 1966), which   implies that
 \begin{equation}
 \sum_{\ell=1}^{q_n}\vert \wtheta_\ell -\theta_\ell \vert \leq
 \frac{q_n}{\sqrt n}\,\sqrt n \| \wbUps_{\weigh}-\bUpsi_{\weigh}\|_{\itG^{k-1}}\;.
 \label{cota}
 \end{equation}
 On the other hand, we have
 $$\esp\left[|\widehat{\itU}_n -{\itU}|\vert \wtX_n\right]=\esp\left[|\widehat{\itU}_n -{\itU}_n+{\itU}_n-{\itU}|\;\vert \tilde
 X_n\right]\leq \sum_{\ell=1}^{q_n}\vert \wtheta_\ell -\theta_\ell \vert +\sum_{\ell>q_n}\theta_\ell $$
 which together with  (\ref{cota}), the fact that $\sqrt n \| \wbUps_{\weigh}-\bUpsi_{\weigh}\|=O_{\prob}(1)$,
 $q_n/\sqrt{n}\to 0$ and  $\sum_{\ell\,\ge 1} \theta_\ell <\infty$ implies that
 \begin{equation}
 \esp\left[|\widehat{\itU}_n -{\itU}|\;\vert \wtX_n \right]\convprob 0\,.
 \label{convergeitU}
 \end{equation}
  We also have the following inequalities
 \begin{eqnarray*}
 \prob(\widehat{\itU}_n\leq t\vert \wtX_n)&=& \prob(\widehat{\itU}_n\leq t \cap |\widehat{\itU}_n-{\itU}|<\epsilon\;\vert \wtX_n)+ \prob(\widehat{\itU}_n\leq t \cap |\widehat{\itU}_n-{\itU}|>\epsilon\;\vert \wtX_n)\\
 &\le & \prob( {\itU}\leq t +\epsilon)+ \prob(|\widehat{\itU}_n-{\itU}|>\epsilon\;\vert \wtX_n)\\
 &\le & F_{\itU}( t +\epsilon)+ \frac 1\epsilon \esp(|\widehat{\itU}_n-{\itU}|\;\vert \wtX_n)\le F_{\itU}(t)+\Delta_\epsilon(t)+\frac 1\epsilon \esp(|\widehat{\itU}_n-{\itU}|\;\vert \wtX_n)\;,
 \end{eqnarray*}
 where $\Delta_\epsilon(t)=\sup_{\vert\delta\vert \leq \epsilon}\vert
 F_{\itU}( t+\delta ) -F_{\itU}(t)\vert$.
  Besides,
 \begin{eqnarray*}
 \prob(\widehat{\itU}_n\leq t\;\vert \wtX_n)&=& \prob(\widehat{\itU}_n\leq t \cap |\widehat{\itU}_n-{\itU}|<\epsilon\;\vert \wtX_n)+ \prob(\widehat{\itU}_n\leq t \cap |\widehat{\itU}_n-{\itU}|>\epsilon\;\vert \wtX_n)\\
 &\ge & \prob( {\itU}\leq t-\epsilon \cap |\widehat{\itU}_n-{\itU}|<\epsilon\;\vert \wtX_n)\\
 &\ge & F_{\itU}( t -\epsilon)- \frac{1}{\epsilon} \esp(|\widehat{\itU}_n-{\itU}|\;\vert \wtX_n)\ge F_{\itU}(t) -\Delta_\epsilon(t)-\frac 1\epsilon \esp(|\widehat{\itU}_n-{\itU}|\;\vert \wtX_n)\;.
 \end{eqnarray*}
 Therefore,
 $$|\prob(\widehat{\itU}_n\leq t\;\vert \wtX_n)-F_{\itU}(t)|\le \Delta_\epsilon(t)+\frac 1\epsilon \esp(|\widehat{\itU}_n-{\itU}|\;\vert \wtX_n)\;.$$
 As we  mentioned in Remark \ref{stest}.1,  $F_{\itU}$ is a continuous
 distribution function on $\real$  and so, uniformly continuous, hence $\lim_{\epsilon\to 0}\,\sup_{t\in \real}
 \;\Delta_\epsilon(t)=0$, which together with (\ref{convergeitU}) implies that
 $\rho_{\mbox{\footnotesize\sc k}}( F_{{\itU}^*_n\vert \wtX_n},
 F_{\itU})=\sup_{t}|\prob(\widehat{\itU}_n\leq t\;\vert \tilde
 X_n)-F_{\itU}(t)|\convprob 0$. \square

 \end{document}